\documentclass[11pt,reqno]{amsart}

\usepackage{cite,epic,eepic,euscript,verbatim,amsmath,amssymb,amsthm,afterpage,float,bm,color} 
\usepackage{graphicx,epsfig}
\usepackage{comment}

\newtheorem{conjecture}{Conjecture}

\numberwithin{equation}{section}

\numberwithin{figure}{section}
\numberwithin{table}{section}

\def\XXint#1#2#3{{\setbox0=\hbox{$#1{#2#3}{\int}$}
\vcenter{\hbox{$#2#3$}}\kern-.51\wd0}}

\def\XXint#1#2#3{{\setbox0=\hbox{$#1{#2#3}{\int}$}
     \vcenter{\hbox{$#2#3$}}\kern-.5\wd0}}

\title[Generalized BO equation]
{Dynamics of solutions in the generalized \\
Benjamin-Ono equation: a numerical study}

\author[S. Roudenko]{Svetlana Roudenko}
\address{Department of Mathematics \& Statistics\\Florida International University,  Miami, FL, USA}
\curraddr{}
\email{sroudenko@fiu.edu}

\author[Z. Wang]{Zhongming Wang}
\address{Department of Mathematics  \& Statistics\\Florida International University,  Miami, FL, USA}
\curraddr{}
\email{zwang6@fiu.edu}

\author[K. Yang]{Kai Yang}
\address{Department of Mathematics \& Statistics\\Florida International University,  Miami, FL, USA}
\curraddr{}
\email{yangk@fiu.edu}



\begin{document}

\begin{abstract}
We consider the generalized Benjamin-Ono (gBO) equation on the real line,
$ u_t + \partial_x (-\mathcal H u_{x} + \tfrac1{m} u^m) = 0, x \in \mathbb R, m = 2,3,4,5$,
and perform numerical study of its solutions. We first compute the ground state solution to $-Q -\mathcal H Q^\prime +\frac1{m} Q^m = 0$ via Petviashvili's iteration method. We then investigate the behavior of solutions in the Benjamin-Ono ($m=2$) equation for initial data with different decay rates and show decoupling of the solution into a soliton and radiation, thus, providing confirmation to the soliton resolution conjecture in that equation. In the mBO equation ($m=3$), which is $L^2$-critical, we investigate solutions close to the ground state mass, and, in particular, we observe the formation of stable blow-up above it. Finally, we focus on the $L^2$-supercritical gBO equation with $m=4,5$. In that case we investigate the global vs finite time existence of solutions, and give numerical confirmation for the dichotomy conjecture, in particular, exhibiting blow-up phenomena in the supercritical setting.
\end{abstract}

\maketitle

\tableofcontents

\section{Introduction}
We study the generalized Benjamin-Ono (gBO) equation
\begin{equation}\label{gBO}
\text{(gBO)} \qquad
u_t + \partial_x (-\mathcal H u_{x} + \tfrac1{m} u^m) = 0, \quad x \in \mathbb R, ~t \in \mathbb R, ~ m \in \mathbb Z^+, 
\end{equation}
where the Hilbert transform $\mathcal H$ is defined by
\begin{equation}\label{HT}
\mathcal H  f(x) = \frac{1}{\pi}\, 
\mbox{p.v.} \int_{-\infty}^\infty \frac{f(y)}{x-y} \, dy,
\end{equation}
and  $\widehat{\mathcal H  f} (\xi) = - i \, \mbox{sgn}(\xi) \, \hat{f}(\xi)$.
We note that $\mathcal H \partial_x = \mathcal D^1$, where $\widehat {\mathcal D^\alpha f}(\xi) = |\xi|^\alpha\, \hat{f} (\xi)$, $\alpha \in \mathbb R$, is the Riesz potential of order $-\alpha$. 
This equation is a generalization of the well-known Benjamin-Ono (BO) equation when $m=2$,\begin{equation}\label{BO}
\text{(BO)} \qquad \quad \qquad
u_t -\mathcal H u_{xx} + u_x u =0, \qquad \qquad \qquad \qquad
\end{equation}
derived by Benjamin \cite{benjamin_1967} and later by Ono \cite{Ono_1975} to model one-dimensional waves in deep water. Other nonlinearities are also relevant in various wave models, e.g., see \cite{davis_acrivos_1967}, \cite{ABFS}, \cite{BK-1}, \cite{BK-2}. The equation \eqref{gBO} with $m=3$ is typically referred to as the {modified} Benjamin-Ono (mBO) equation.

The gBO equation formally conserves several quantities: the $L^2$ norm (often called mass), the energy (or Hamiltonian), and the $L^1$-type integral
\begin{align}
\label{E:M}
& M[u(t)]  = \int u(t)^2 \, dx = M[u(0)], \\
\label{E}
& E[u(t)] =\dfrac{1}{2}\int |\mathcal D^{1/2} u(t)|^2 - \dfrac{1}{m(m+1)}\int u(t)^{m+1} = E[u(0)], \\
\label{E:L1}
& \int u(x,t) \, dx  = \int u(x,0) \, dx.
\end{align}
We note that in the gBO hierarchy, the BO equation \eqref{BO} is completely integrable, for example, the Lax pair was originally constructed in \cite{N_79} and \cite{BK_79}, see also \cite{AF1983}, \cite{KM_98}.

Besides other symmetries, the gBO equation is invariant under the scaling:
if $u(x,t)$ is a solution to \eqref{gBO}, then so is
\begin{equation}\label{E:scaling}
u_\lambda(x,t)=\lambda^{\frac{1}{m-1}} u(\lambda x, \lambda^2 t), \, \lambda > 0.
\end{equation}
Under this symmetry the Sobolev $\dot{H}^s$ norm with $s = \frac12- \frac1{m-1}$ is invariant.
The value of this scaling critical index $s$ classifies the equation \eqref{gBO} as follows: when $m=2$ ($s<0$), the equation \eqref{BO} is called the $L^2$-subcritical equation; if $m=3$ ($s=0$), the equation is $L^2$-critical;  when $m \geq 4$ ($s>0$), the equation \eqref{gBO} is $L^2$-supercritical. This classification is useful for understanding the long-term behavior of solutions to \eqref{gBO} with various power nonlinearities, which is the goal of this work.

The wellposedness theory for the Cauchy problem for the BO equation \eqref{BO}
was initiated by Saut in \cite{S1979}, followed by \cite{I1986},  \cite{ABFS}, \cite{GV1989}, \cite{P1991}, with further improvements in \cite{KT2003}, \cite{KK2003}, \cite{T2004}, \cite{BP2008}, finally resulting in the global well-posedness in the $L^2$ space by Ionescu \& Kenig in \cite{IK2007} (see also \cite{MP2012} and \cite{IT2019}).
The well-posedness theory for the generalized equation \eqref{gBO}, for example, goes back to the work of Kenig, Ponce \& Vega \cite{KPV1994}, further improvements were done in \cite{MR2004a}, \cite{MR2004b}, \cite{BP2006}, \cite{V2009}; the best available result for the mBO ($m=3$) is in the energy-space $H^{1/2}$ by Kenig \& Takaoka \cite{KT2006}, and by Vento \cite{V2010} and Molinet \& Ribaud \cite{MR2004b} for gBO with $m=4$ in $H^s$ for $s >\frac13$ and with $m\geq 5$ for $s \geq s_m = \frac12-\frac1{m-1}$.  For the purpose of this paper, it suffices to have local well-posedness in the energy space $H^{1/2}$.
From the local theory it follows that solutions to the gBO equation \eqref{gBO} have a maximal forward lifespan $[0,T)$ with either $T=+\infty$ or $T<+\infty$ (in the energy space or in the space as discussed above). If $T < + \infty$, then $\lim_{t \nearrow T} \| \mathcal D^{1/2} u (t) \|_{L^2} = +\infty$. Furthermore, if the time $T$ is finite, the local well-posedness and scaling provide the lower bound $\|\mathcal D^{1/2} u(t)\|_{L^2} \gtrsim (T-t)^{-\frac1{2(m-1)}}$ for $t$ close to $T$ (e.g., see \cite[Remark 2.4]{MP2016}).

The gBO equation has a family of localized traveling solitary waves, which propagate in the positive $x$-direction
\begin{equation}\label{E:Q}
u(x,t) = Q_c(x-ct+C), ~~ c>0, ~~ C \in \mathbb R,
\end{equation}
where $Q_c = c \, Q(cx)$ and $Q$ is the unique positive, even, decreasing (for $x>0$) and vanishing at infinity solution of
\begin{equation}\label{Q}
- Q - \mathcal H Q^\prime + \frac1{m}Q^{m}=0.
\end{equation}
The existence and uniqueness of this solution in the BO equation ($m=2$) was obtained in \cite{AT1991}, in that case it is explicit $Q(x) = \frac4{1+x^2}$. For a general nonlinear power, the existence of $Q$ follows from \cite{W1987} and \cite{ABS1997}, though for $m \neq 2$, $Q$ is not explicit but still has a polynomial decay $|x|^{-2}$ (e.g., see \cite{BX12}). The uniqueness was established in \cite{FL2013} (in a more general case than \eqref{Q}). To follow the literature, we refer to this unique solution (different for each $m$) as the {\it ground state}.  The stationary problem \eqref{Q} is related to the best constant in the following Gagliardo-Nirenberg inequality
\begin{equation}
\label{E:GN}
\| v \|_{L^{m+1}}^{m+1} \leq C_m \, \|\mathcal D^{1/2} v\|_{L^2}^{m-1} \, \|v\|^{L^2}_2
\end{equation}
with $\displaystyle C_m = \frac{m(m+1)}{2} \left( \frac2{(m-1) \|Q\|_{L^2}^2} \right)^{\frac{(m-1)}2}$, see \cite{ABLS2002} and also \cite{FLP2014}.

In the $L^2$-subcritical case ($m=2$, BO equation) the family of solitary waves is explicit
\begin{equation}\label{soliton}
u(t,x)=\frac{4c}{1+c^2(x-ct+C)^2} \equiv Q_c (x-ct +C),
\end{equation}
where $C \in \mathbb R$ is a location or shift of the ground state $Q$ and $c>0$ is the speed of the soliton and its scaling parameter. 
Note that there can be large amplitude and arbitrarily small amplitude solitary waves. As it was mentioned, the BO equation is integrable, and thus, has an inverse scattering formalism,
e.g., see \cite{AFA_83,BK_79,KM_98,AF1983,KLM_99,KLM_98,Matsuno_04,MX_11,MX_12, Wu2017}.

Since the BO equation \eqref{BO} is subcritical, the solitary waves are stable.
For example, the orbital stability in the energy space $H^{1/2}$ has been shown in \cite{BBSSB}, \cite{BSS}, \cite{ABH}. A more delicate asymptotic stability in $H^{1/2}$ has been obtained in \cite{KM2009} and also in \cite{GTT2009}. In the mBO equation \eqref{gBO}, $m=3$, there are also traveling solitary waves, however, as the equation is $L^2$-critical, the behavior will depend on the size (and shape) of the initial data. Similar to the $L^2$-critical nonlinear Schr\"odinger (NLS) and general KdV (gKdV) equations, any $H^{1/2}$ initial datum $u_0$ with the mass below the soliton mass, $\|u_0\|_{L^2} < \|Q\|_{L^2}$, produces a solution $u(x,t)$ that exists globally with the bounded $H^{1/2}$ norm, which follows from the Gagliardo-Nirenberg inequality \eqref{E:GN} and its sharp constant. At the mass threshold $\|u_0\|_{L^2} = \|Q\|_{L^2}$, solutions may blow-up, those would be called {\it minimal mass} blow-up solutions, and first such example was shown in \cite{MP2016}. This, of course, implies instability of solitary waves in the mBO equation. The minimal mass blow up dynamics is unstable, since a perturbation of it can lead to a global solution (if the initial data less than the mass of $Q$). The stable blow-up is unknown and in this work we show numerical simulations toward that. In other gBO equations $m>3$, or in the $L^2$-supercritical case, the dynamics of solutions and solitary waves is even less understood, though the work of Farah, Linares \& Pastor \cite{FLP2014} establishes the global well-posedness below the `mass-energy' threshold, similar to the NLS case, introduced by the first author and Holmer in \cite{HR2007}, \cite{HR2008}, following the ideas of \cite{KM2006}.
Some numerical investigations for the BO-type equations were done in \cite{KS2015}, see also numerical studies of solitons and interactions in the integrable BO equation in \cite{PD00} and \cite{MPST93}. For further details, see excellent reviews \cite{P2015} and \cite{S2018}.

In this work we start numerical investigations of the gBO equations on the real line to obtain a better understanding of solutions behavior in the $L^2$-critical and $L^2$-supercritical cases, in particular, we study the behavior of solitary waves, evolution of different types of initial data, global existence and formation of stable blow-up. In particular, we investigate the following conjectures.

\begin{conjecture}[$L^2$-critical case, mBO]\label{C: dichotomy}
Let $u_0 \in H^{1/2}(\mathbb R)$ be sufficiently localized and let $Q$ be the ground state solution of \eqref{Q} with $m=3$.
Let $u(x,t)$ be the mBO time evolution of $u_0(x)$.  Then
\begin{itemize}
\item[I.] $u(x,t)$ exists globally in time
if $\|u_0\|_{L^2}<\|Q\|_{L^2}$,
\item[II.] $u(x,t)$ blows up in finite time if $E[u_0]<0$ (which implies $\|u_0\|_{L^2} > \|Q\|_{L^2}$).
\end{itemize}
\end{conjecture}
{\bf Remarks.}
1. The global existence in Part I holds by the argument of Weinstein \cite{W1987}, see for example, \cite{MP2016}.
2. There are also positive energy blow-up solutions (with $\|u_0\|_{L^2} > \|Q\|_{L^2}$).

\begin{conjecture}[$L^2$-supercritical case, gBO, $m > 3$]\label{C: dichotomy2}
Let $u_0 \in X \cap H^{1/2}(\mathbb R)$ be sufficiently localized (here, $X$ is the corresponding local well-posedness space for a given $m$) and let $Q$ be the ground state solution of \eqref{Q} for a given $m$.
Denote $D[u]=\| (\mathcal H \partial_x)^{1/2} u \|_{L^2}^2$ and let $u(x,t)$ be the gBO time evolution of $u_0(x)$. Then
\begin{itemize}
\item[I.]
$u(x,t)$ blows up in finite time if $E[u_0]<0$,

\item[II.]
the following dichotomy holds (recall $s = \frac12 - \frac1{m-1}$):\\
Suppose
\begin{align}\label{ME bound}
M[u_0]^{\frac12-s} E[u_0]^{s} < M[Q]^{\frac12-s} E[Q]^{s}, ~E[u_0]>0.
\end{align}
Then 
\begin{itemize}
\item[(IIa)]
$u(x,t)$ exists globally in time if
\begin{align}\label{MD scatter}
M[u_0]^{\frac12-s} D[u_0]^{s} < M[Q]^{\frac12-s} D[Q]^{s},
\end{align}
\item[(IIb)]
$u(x,t)$ blows up in finite time if
\begin{align}\label{MD blowup}
M[u_0]^{\frac12-s} D[u_0]^{s} > M[Q]^{\frac12-s} D[Q]^{s}.
\end{align}
\end{itemize}
\end{itemize}
\end{conjecture}

{\bf Remark.} The global existence in Part IIa of Conjecture \ref{C: dichotomy2} (or maximal time interval existence in case $m=4,5$), i.e., under the conditions \eqref{ME bound} and \eqref{MD scatter}, was proved by Farah, Linares \& Pastor in \cite{FLP2014} (even in a more general case of dispersion).

In this paper we give positive confirmation of the above conjectures via numerical simulations, where the `sufficiently localized' data has a single maximum and monotone decay from it (we are able to investigate exponential decay and polynomial decay as slow as $1/|x|$). Besides studying the $L^2$-critical and supercritical cases, we also consider the $L^2$-subcritical case ($m=2$). Having an explicit soliton $Q =\frac4{1+x^2}$ in the BO equation is useful for validating numerical approaches and accuracy. Being completely integrable, the BO equation can be studied by the inverse scattering methods, which could provide asymptotics for a smooth class of rapidly decaying solutions, and that can be helpful in understanding the time evolution that is observed in numerical simulations. We run simulations of the BO time evolution for different sizes and types of initial data and observe the following behavior of time evolution of a single peak data:

- if the maximum of the solution starts moving to the right, the solution then converges to the (shifted) soliton $Q_c$ traveling right with the speed approaching $c$ (and some radiation to the left),

- if the maximum of the solution starts moving to the left, in some cases\footnote{We could not check all cases due to the computational restrictions.} it slows down, then reverses the direction and  travels to the right, emerging the soliton; the time evolution will decouple into the soliton component approaching the rescaled (and shifted) soliton $Q_c$ traveling right with the speed $c$  plus radiation moving to the left.

{\bf Remark.} The behavior in the second part shows that the solution while starting to radiate, eventually decouples into a (rescaled) soliton traveling to the right and radiation going to the left. Due to numerical constrains we could not so far run the simulations long enough for all considered types of data (as in \eqref{E:dataset}) to check if the solution always `finds' a soliton, which will then decouple from radiation.
Of course, such behavior is expected not only in the completely integrable setting  but for general dispersive PDEs (e.g., soliton resolution conjecture \cite{T2004b} or grand conjecture \cite{S2006}). In the BO equation it is theoretically possible to solve the Cauchy problem by a direct scattering for a smooth rapidly decaying initial data (see \cite{Wu2017}). Here, we are able not only to investigate the asymptotic behavior, but also observe the process of soliton formation or decoupling of the soliton and the radiation; furthermore, our numerical methods (of rational basis function) allowed us to investigate slow decaying initial conditions (such as ~$1/|x|$).

We now mention the numerical methods of this paper.
One of the main numerical difficulty lies in the evaluation of the Hilbert transform $\mathcal H$. Due to its nature, {\it spectral} expansions are often used in numerical studies, for which fast Fourier transform (FFT) may be applied. In the literature, the following methods exist: Fourier basis method \cite{DM09,Weideman95,JW92, MPST93,DDM16}, radial basis function method \cite{Fa07} and rational basis function method \cite{Weideman95,JW92,BX12}. A numerical comparisons of difference types of these expansions is given in \cite{BX11}. Other types of discretization for the BO equation, such as finite difference  \cite{TM98,DHKR16}  and Galerkin finite element \cite{So18}, have also been introduced and applied. Temporal discretization methods, such as leapfrog \cite{Weideman95} and Crank-Nicolson \cite{So18}, are used for simulating the dynamics of solutions to the BO and gBO equations via standard explicit fourth order Runge-Kutta (RK) method, such as RK4. To circumvent the strict CFL (Courant-Friedrichs-Lewy) conditions of explicit methods, the exponential time differencing (ETD) \cite{BKJV98} was developed. We mention that the ETD schemes coupled with RK4 are proposed in \cite{CM02} for stiff systems (referred as ETDRK4), and further improved in \cite{KT05} for stability (which are referred as mETDRK4).

The goal of this project is to numerically study the BO \eqref{BO}, mBO and gBO \eqref{gBO}  equations and investigate solitary waves, scattering and blow-up phenomena on the {\it whole real line}. The designed numerical schemes are first verified by the explicit solitary waves solutions for \eqref{Q} and \eqref{BO}, and then employed to simulate solution dynamics such as scattering, solitary waves, or blow-ups. In order to accurately capture the non-solitary waves, we use both Fourier basis and rational \eqref{BO} basis spectral expansions in space, and  fourth order explicit schemes in time such as the standard fourth order Runge-Kutta (RK4) method. During long time simulations, a modified exponential time differencing (ETD) coupled with RK4 is also used for efficiency by relaxing the strict CFL restriction in RK4.

This paper is organized as follows: in Section \ref{S:spatial} we describe the spatial discretizations used in our numerical algorithms, in Section \ref{S:groundstate} we show how we use the discretizations and the Petviashvili's iterations to obtain the ground state solutions to \eqref{Q} for $m=2,3,4,5$. In Section \ref{S:dynamics} we first discuss the temporal discretization in order to simulate solutions to the dynamical equation with various checks of consistency, reliability, and errors, and then perform numerical simulations and analysis: for the BO equation in \S \ref{S:BO}, giving confirmation to the soliton resolution conjecture in that equation, for the mBO equation in \S \ref{S:mBO}, giving numerical confirmation of Conjecture \ref{C: dichotomy}, and finally in \S \ref{S:gBO} we study the gBO equation, $m=4,5$, and give numerical confirmation for Conjecture \ref{C: dichotomy2}.

{\bf Acknowledgments.}
S.R. was partially supported by the NSF grant DMS-1927258 as well as part of the K.Y.'s research and travel support on this project came from the same grant.

\section{Spatial discretization} \label{S:spatial}
One of the main difficulties in numerical simulations of the BO or gBO equations  is in dealing with the {\it infinite} domains of the equations. There are mainly two types of treatments, see \cite{JW92}. The first one is to restrict the equations \eqref{BO} or \eqref{gBO} onto a large closed domain $[-L, L]$ with suitable boundary conditions. This setting is viable if waves are not present at far field. The second one is to use a change of variable $x=L\tan\frac{\theta}{2}$ for the improper integration over the infinite domain.  Spectral methods are typically preferred in both treatments due to their high order accuracy and small values of the solution near the  boundary.

\subsection{Fourier expansion}
When the equations \eqref{BO} and \eqref{gBO} are mapped onto a closed domain $[-L, L]$, it is numerically simple and efficient to consider the periodic boundary conditions, for which the FFT may be applied. The periodic boundary conditions may be applied for solitary wave type solutions for large enough size of the domain $L$.

With a uniform partition of $N$ grid points, $x_j=jh$, where $h=L/N, j=-N,\cdots, N$, and the standard discrete Fourier expansion, one obtains
\begin{equation} \label{Fourier}
u_j:=u_j(t,x_j)=\sum_{n=-N}^{N-1}a_ne^{i n\pi x_j/L}, \text{  where } a_n=\frac{1}{2N}\sum_{j=-N}^{N-1} u_j e^{-in\pi x_j/L}.
\end{equation}
By introducing vectors $u=(u_{-N},\cdots, u_{N-1})^T, \quad a=(a_{-N},\cdots, a_{N-1})^T$ and matrices
$$
F^{-1}_{jn}=e^{in\pi x_j/L},\quad F_{nj}=\frac{1}{2N}e^{-in\pi x_j/L}, \quad -N\leq j,n\leq N-1, $$
One can write $u$ in a matrix form as $u=F^{-1}a, \quad a=Fu$.

Using properties of Fourier transform, one can easily obtain
\begin{align}
u'_j &= \frac{i\pi}{L} \sum_{n=-N}^{N-1}n a_ne^{i n\pi x_j/L},\\
u''_j &=-\frac{\pi^2}{L^2} \sum_{n=-N}^{N-1}n^2 a_ne^{i n\pi x_j/L},\\
\mathcal H(u'')_j &=-\frac{i\pi^2}{L^2}\, \text{sgn}(n) \sum_{n=-N}^{N-1}n a_ne^{i n\pi x_j/L},
\end{align}
where $\text{sgn}(n)$ is the standard sign function.
The matrix forms of the derivatives and the Hilbert transform are
$$
u'=F^{-1}E_1Fu,\quad u''=F^{-1}E_2Fu,\quad \mathcal H(u)=F^{-1}E_{\mathcal H} Fu,
$$
where
\begin{align}
E_1&=\frac{i\pi}{L} \, \text{diag}(-N,\cdots, (N-1)), \\
E_2 &=-\frac{\pi^2}{L^2} \, \text{diag}(N^2,\cdots, (N-1)^2,\cdots, 0, \cdots, (N-1)^2),\\
E_{\mathcal H} &=-\frac{\pi^2}{L^2} \, \text{diag}(-1,\cdots, -1, 1,\cdots, 1).
\end{align}
We note that the matrices $F$ and $F^{-1}$ are usually not explicitly assembled, instead an FFT algorithm is typically used in computation.

The semi-discrete form of the BO and gBO equations is
\begin{equation}\label{semidiscreteFm}
u_t+ F^{-1}E_1F B(u) +F^{-1}E_{\mathcal H} E_2 F u=0,
\end{equation}
where the  matrix function
$B(u)=\frac{1}{m} \text{diag}(u^m_{-N}, \cdots, u^m_{N-1})$ is used for the BO \eqref{BO} with $m=2$, and gBO \eqref{gBO} for other $m$.

\subsection{Rational basis expansion}
The second approach to treat an infinite domain is to make a change of variable $x=L\tan\frac{\theta}{2}$.  We  first consider the following rational basis expansion
\begin{equation}\label{rational}
 u(t,x)=\sum_{n=-\infty}^{\infty} a_n(t)\phi_n(x), \quad \phi_n(x)=\frac{(L+ix)^n}{(L-ix)^{n+1}},
\end{equation}
where $L$ is a mapping parameter indicating that half of the grid points are located in the interval $[-L,L]$. This parameter is to be chosen depending on accuracy expectation. A simple computation shows that  $\{ \phi_n(x) \}_{n=-\infty}^{\infty}$ 
forms a complete orthogonal basis in $L^2(-\infty, \infty)$ with the following orthogonality
normalization (see also \cite{Christov82})
$$
\int_\infty^\infty \phi_m(x) \, \overline{\phi_n(x)} \, dx = \begin{cases}
      \pi/L, &  m=n,\\
      0, & m\neq n.
    \end{cases}
$$
Therefore, we have
$$
a_n(t)=\frac{L}{\pi}\int_{n=-\infty}^{\infty} u(t,x) \, \phi_n(x) \, dx.
$$
Note that the solitary waves $Q_c$ 
can be expressed by two modes in the rational expansion:
$$
\frac{4c}{1+c^2x^2}= 2 \, \phi_{-1}(x)+ 2 \, \phi_0(x)
$$
with $L=1/c$.
Furthermore, using the rational expansion \eqref{rational}, the Hilbert transform \eqref{HT} can be easily calculated, see \cite{Weideman95}, as
$$
\mathcal H(u(t,x))=\sum_{n=-\infty}^{\infty} i \, a_n(t) \, \text{sgn}(n) \, \phi_n(x).
$$
The derivatives of $u(t,x)$ can also be easily computed as 
\begin{align}
u_x(t,x)=\sum_{n=-\infty}^{\infty} &\frac{i}{2L} \, \big[ na_{n-1}+(2n+1)a_n+(n+1)a_{n+1} \big] \, \phi_n(x), \label{ux}\\
u_{xx}(t,x)=\sum_{n=-\infty}^{\infty} & -\frac{1}{4L^2} \, \big[n(n-1)a_{n-2}+4n^2a_{n-1}+(6n^2+6n+2) a_n  \notag\\
             & +4(n+1)^2a_{n+1} +(n+2)(n+1)a_{n+2} \big] \, \phi_n(x). \label{uxx}
\end{align}
For the numerical computations, a truncation of $2N$ terms is used, i.e.,
$$
u(t,x) \cong a^T\phi:=\sum_{n=-N}^{N-1} a_n(t) \, \phi_n(x),
$$
which leads to the matrix forms
$$
\mathcal H(u)=[C_{\mathcal H} a]^T\phi, \quad u_x=[C_1a]^T\phi, \quad u_{xx}=[C_2a]^T\phi,
$$
where $C_{\mathcal H} = i \, \text{diag}(-1, \cdots, -1, 1, \cdots, 1)$, and $C_{1,2}$ are given in \eqref{ux} and \eqref{uxx}.

Now by a change of variable $ x=L\tan \frac{\theta}{2}, ~~ -\pi\leq \theta \leq \pi,$
and the spatial discretization $x_j=L\tan \frac{\theta_j}{2}$, where $\theta_j=jh, h=\pi/N, j=-N, \cdots, N$, we obtain
\begin{equation}\label{uj}
u_j=u(t,L-ix_j)=\sum_{n=-N}^{N-1}a_n e^{i\theta_j}.
\end{equation}
We note that the above discretization in space is not uniform in $x$, but uniform in $\theta$. Also, the truncation is used only in the spectral expansion and not in space.

Equipped with \eqref{uj}, the semi-discrete version of the BO equation \eqref{BO} or gBO equation \eqref{gBO} is
\begin{equation}\label{semidiscreteRa}
u_t+ P^{-1}F^{-1}C_1FPB(u) +P^{-1}F^{-1}C_{\mathcal H} C_2FPu=0,
\end{equation}
where  $F$ is the standard FFT matrix with $\{ n\theta_j \}$, and $P=\text{diag}(L-ix_{-N}, \cdots, L-ix_{N-1})$ is the weight matrix, which comes from the rational basis.
Now that we described approaches for spatial discretization, we will discuss how we obtain the ground state solution to the nonlinear equation \eqref{Q}.

\section{Ground state solutions} \label{S:groundstate}
\subsection{Iterative scheme}

To obtain the ground state solution $Q$ of \eqref{Q}, we use the Petviashvili's iteration  method from \cite{PS2004}, which is also called the {\it renormalization method}. The convergence analysis for this scheme can be found in \cite{PS2004}, \cite{OSSS2016}.

Taking \eqref{Q} and multiplying it 
by $Q$, then integrating over the real line $\mathbb{R}$, one obtains the two quantities $SL(Q)$ and $SR(Q)$ (on the left and right sides of the equation)
\begin{align}\label{GBO invariant}
SL(Q):=\int_{\mathbb{R}} Q^2 dx=\int_{\mathbb{R}} Q \, \left[(\mathcal H \partial_x +1)^{-1} \frac{Q^m}{m} \right] \, dx=: SR(Q).
\end{align}
To prevent the fixed point iteration scheme
\begin{equation}\label{E:Qscheme}
Q^{(n+1)}=(\mathcal H \partial_x +1)^{-1} \dfrac{ \left({Q^{(n)}} \right)^m}{m}
\end{equation}
from going to $0$ or $\infty$, we search for a constant $c_n$ at each iteration such that the fixed point iteration $c_nQ^{(n)}$ preserves the equality \eqref{GBO invariant}, i.e.,
\begin{align}\label{GBO invariant c}
SL(c_n \, Q^{(n)})=SR \left(c_n \, {Q^{(n)}} \right).
\end{align}
Inserting \eqref{GBO invariant c} into \eqref{GBO invariant} yields
\begin{align}\label{GBO invariant c2}
c_n= \left( \dfrac{SL(Q^{(n)})}{SR(Q^{(n)})} \right)^{\frac{1}{m-1}}.
\end{align}
Therefore, we have the iteration
\begin{align*}
Q^{(n+1)}=(\mathcal H\partial_x +1)^{-1} \dfrac{ \left(c_n \, {Q^{(n)}} \right)^m}{m}.
\end{align*}
Putting \eqref{GBO invariant c2} in the above equation leads to the iterative numerical scheme
\begin{align}\label{GBO iteration}
Q^{(n+1)}=  \dfrac{1}{m} \left( \dfrac{SL(Q^{(n)})}{SR(Q^{(n)})} \right)^{\frac{m}{m-1}} (\mathcal H \partial_x +1)^{-1} \left({Q^{(n)}} \right)^m.
\end{align}
Therefore, using different discretizations described in the previous section, the numerical iterative schemes become
\begin{align}\label{GBO Q Fourier}
Q^{(n+1)}=  \left( \dfrac{SL(Q^{(n)})}{SR(Q^{(n)})} \right)^{\frac{m}{m-1}} F^{-1} \left(E_{\mathcal H} E_1+I \right)^{-1}F B(Q^{(n)})
\end{align}
for the Fourier discretization, and
\begin{align}\label{GBO Q rational}
Q^{(n+1)}=  \left( \dfrac{SL(Q^{(n)})}{SR(Q^{(n)})} \right)^{\frac{m}{m-1}} P^{-1}F^{-1} \left(C_{\mathcal H} C_1+I \right)^{-1}FP B(Q^{(n)})
\end{align}
for the rational basis expansion discretization. Here $I$ is the identity matrix and the matrix
$B(Q) = \frac{1}{m} ~\mbox{diag} ~~ (Q^m_{-N}, \cdots, Q^m_{N-1})$.

We stop our fixed point iteration \eqref{GBO iteration} when $\|Q^{(n+1)}-Q^{(n)}\|_{\infty}<Tol$,  with a typical  tolerance $Tol=10^{-12}$.

\subsection{Profiles of ground states}\label{3.2.1}
Here, we compute the ground state solutions $Q$ of \eqref{Q} for $m=2,3,4,5$ by using the iterative method discussed above. Unless specified otherwise, we employ the rational basis functions with the number of nodes $N=4096$, the mapping parameter $L=20$, and the initial guess $Q^{(0)}=e^{-x^2}$.

\begin{figure}[ht]
\includegraphics[width=0.49\textwidth]{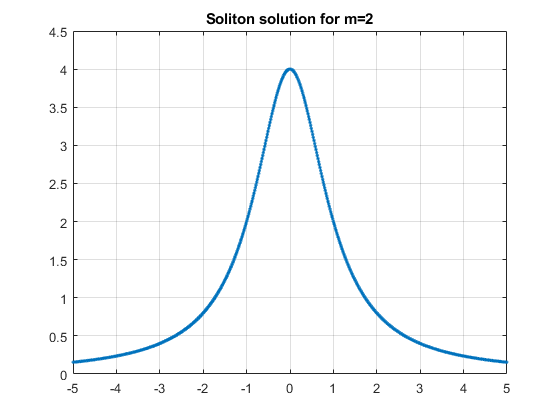}
\includegraphics[width=0.49\textwidth]{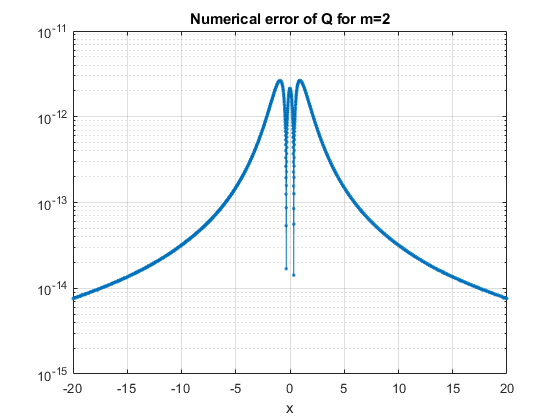}
\caption{The ground state $Q$ in the BO ($m=2$) equation. Left: the computed $Q^{num}$ and the exact ground state $Q = \frac4{1+x^2}$ are indistinguishable. Right: the numerical error
$|Q^{num}-Q|$, note that the error is on the order of $\sim 10^{-12}$ or less.}
\label{Q profile}
\end{figure}

In Figure \ref{Q profile} we show our computations for the BO ground state $Q$, $m=2$. On the left we plot both the numerically obtained ground state $Q^{num}$ and the exact solution $Q = \frac4{1+x^2}$, which completely coincide, and thus, indistinguishable. To further confirm it, on the right of Figure \ref{Q profile}, we show the difference between the numerical solution and the exact one, noting that it is on the order of $10^{-12}$ or smaller, which indicates that the numerical computation of $Q$ matches the exact solitary wave $Q=\frac{4}{1+x^2}$ almost up to the machine precision.

\begin{figure}[ht]
\includegraphics[width=0.8\textwidth, height=0.3\textheight]{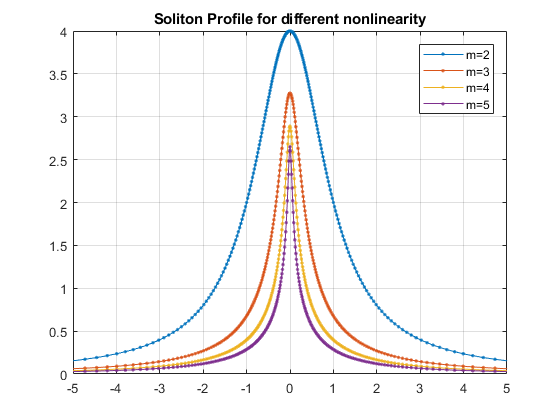}
\caption{The ground state solutions of \eqref{Q} for the gBO equation with $m=2, 3, 4, 5$.}
\label{QprofileCompare}
\end{figure}

Next we compute the ground state solution of \eqref{Q} for other values of $m$. In Figure \ref{QprofileCompare}, we plot the computed ground state profiles for $m = 2,3,4,5$ to provide a comparison for different nonlinear powers. In the case $m=5$, since the solution is more concentrated at the origin, we take $L=10$. Observe that with higher nonlinear power $m$, the height of the ground state decreases, the ground state becomes less flat at the maximum, and gets more and more localized.

We remark that compared to the methods used in \cite{BX12}, the  iterative method allows a more robust initial guess, for example, starting from a bump function such as Gaussian $Q^{(0)}=e^{-x^2}$, the computations will always converge to the ground state solution $Q$ in all considered cases $m=2,3,4,5$. We also used the Fourier spectral discretization and obtained similar results.

\subsection{Pohozaev identities and consistency check} Another way to check the validity of computed ground states is to use Pohozaev identities
\begin{align}\label{Q identity}
\dfrac{2}{m(m+1)} \|Q\|_{L^{m+1}}^{m+1}=\|Q\|_{L^2}^2 \quad \mbox{and} \quad 
\| (\mathcal H \partial_x)^{\frac{1}{2}} \, Q \|_{L^2}^2 = \dfrac{m-1}{2} \| Q \|_{L^2}^2.
\end{align}
They are derived in a similar manner as in the nonlinear Schr\"odinger (NLS) equation  by multiplying the equation \eqref{Q} by $Q$ or $x Q^\prime$. 

Denoting $Q^{num}$ the numerically computed ground state solution from the scheme \eqref{GBO iteration}  and defining the values
\begin{align}
\label{E:e1}
e_1 & = \big| \|Q^{num}\|_{L^2}^2 - \frac{2}{m(m+1)} \|Q^{num}\|_{L^{m+1}}^{m+1} \big|, \\
e_2 & = \big| \frac{m-1}{2} \| Q^{num} \|_{L^2}^2 - \| (\mathcal H \partial_x)^{\frac{1}{2}} Q^{num} \|_{L^2}^2 \big|
\label{E:e2}
\end{align}
as the errors appropriate for the consistency check and validity of the ground state computation, we show both errors in Table \ref{Table: Q error}
for the numerical ground state solution $Q^{num}$ obtained in the previous section \S\ref{3.2.1}. We observe that the quantities are consistent almost to the machine error ($\sim10^{-11}$) in all computed cases of $m$, thus, also re-assuring the later computations of the conserved quantities such as mass and energy.

\begin{table}[ht]
\begin{tabular}{|l|l|l|l|l|}
\hline
            &  $m=2$ & $m=3$ & $m=4$ & $m=5$  \\
\hline
 $e_1$      & $7.5211e-12$ & $1.7453e-12$ & $6.9766e-13$ & $8.8152e-14$  \\
 \hline
 $e_2$      & $1.6707e-11$ & $8.4839e-12$ & $8.4350e-12$ & $5.9854e-12$  \\
 \hline
\end{tabular}
\caption{Consistency check of $Q^{num}$ via Pohozaev identities \eqref{E:e1}-\eqref{E:e2}.} 
\label{Table: Q error}
\end{table}

\section{Dynamical equations}\label{S:dynamics}
After obtaining solutions to the stationary equations and discussing the spatial discretizations used, we now discuss our numerical approaches to track the time evolution in gBO equations and the results that we obtained in our simulations.

\subsection{High order temporal discretizations}
In order to match the high order of accuracy of spatial discretization by spectral methods, high order temporal discretization should be used. For example, both the semi-discrete equations \eqref{semidiscreteFm}  and \eqref{semidiscreteRa} can be solved by a standard fourth order Runge-Kutta method (RK4) for $u_t=f(u)$
\begin{equation}\label{RK4}
u^{k+1}  =  u^{k}+\tfrac{1}{6}(a+b+c+d),
\end{equation}
$$
a  =  \Delta t f(u^k), ~
b  =  \Delta t f(u^k  + a/2), ~
c  =  \Delta t f(u^k+b/2), ~
d  =  \Delta t f(u^k+c),
$$
where $u^k=u(t_k,\cdot)$, $t_k=k\Delta t$, and $\Delta t$ is the temporal step size. We use the fourth order accuracy in the time discretization scheme, same order of accuracy as in our spatial discretizations in \S \ref{S:spatial}.

We note that due to its explicit discretization, the RK4 suffers from a strict CFL
condition, $\Delta t= \mathcal{O}(1/N^2)$, which could be computationally expensive in long time simulations. One solution to this problem for the Fourier method \eqref{semidiscreteFm} is to use the mETDRK4 \cite{KT05}, which is still stable when $\Delta t$ is taken to be roughly 100 times of those in RK4. Note that mETDRK4 is applicable in \eqref{semidiscreteFm} since all $E_1$, $E_2$ and $E_{\mathcal H}$ are diagonal matrices. However, since $C_1$ and $C_2$ in the rational basis expansion method \eqref{semidiscreteRa} are tridiagonal and pentadiagonal, respectively, mETDRK4 is very inefficient, as the exponential of the sparse matrices will result in the full matrices. Therefore, unless otherwise specified, we use mETDRK4 for the Fourier method \eqref{semidiscreteFm} and the standard RK4 method for the rational base expansion method \eqref{semidiscreteRa} in our simulations.

\subsubsection{Consistency check}

\begin{figure}[ht]
\includegraphics[width=0.49\textwidth]{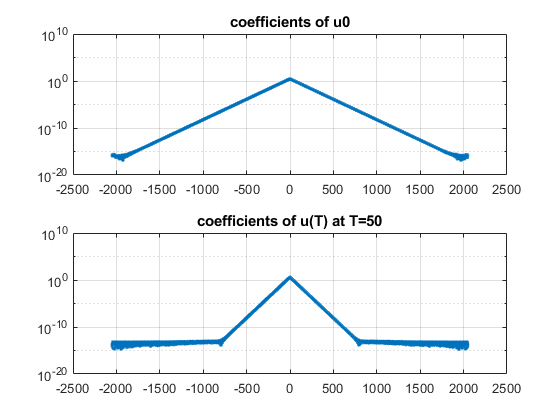}
\includegraphics[width=0.49\textwidth]{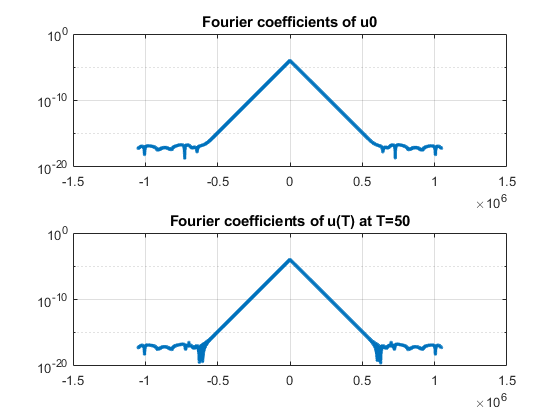}
\caption{Case $m=2$ (BO) with the soliton initial condition $u_0=\frac{4}{1+(x+25)^2}$. Left: coefficients $a_n(t)$ in the rational basis expansion \eqref{rational}. Right: coefficients $a_n(t)$ in the Fourier basis expansion \eqref{Fourier}. Top plots at time $t=0$ and bottom plots at time $t=50$.}
\label{Q coefficients}
\includegraphics[width=0.95\textwidth, height = 0.15\textheight]{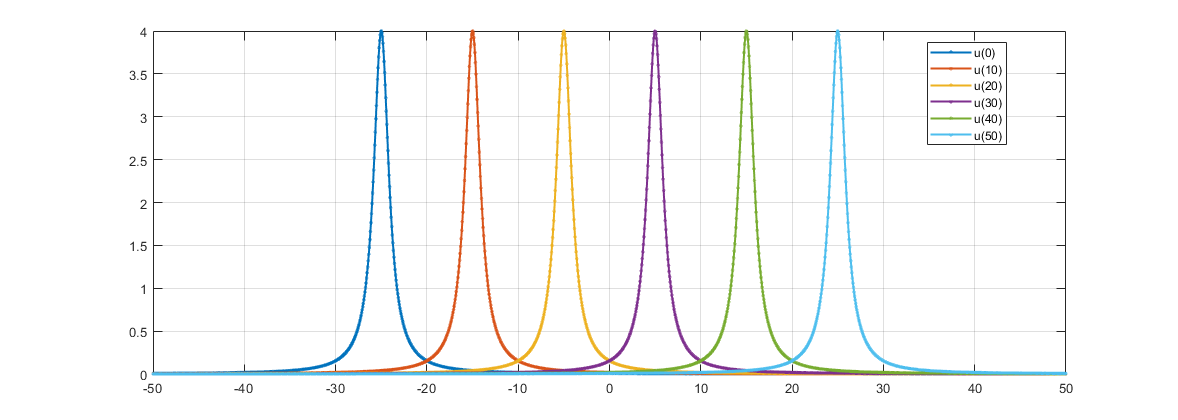}
\caption{Case $m=2$ (BO), the soliton propagation with $u_0=\frac{4}{1+(x+25)^2}$
at $t=0,10,20,30,40,50$.}
\label{Q soliton}
\end{figure}

We take $m=2$ and consider the BO equation \eqref{BO} and one of its soliton solutions $u_0=\frac{4}{({1+(x-x_0)^2})}$. By tracking the coefficients $a_n(t)$ in time, one can observe  the two types of spatial discretization, \eqref{semidiscreteFm} and  \eqref{semidiscreteRa},  are almost equally efficient. In Figure \ref{Q coefficients} we observe that the coefficients of both Fourier and rational basis expansions decay very fast when the number of terms $N$ increases even at large times. This also shows that our choice of $N=4096$ is sufficient in most of our simulations. In Figure \ref{Q soliton} we show the snapshots of the soliton moving to the right in the time evolution of the BO equation.

\begin{figure}[ht]
\includegraphics[width=0.49\textwidth]{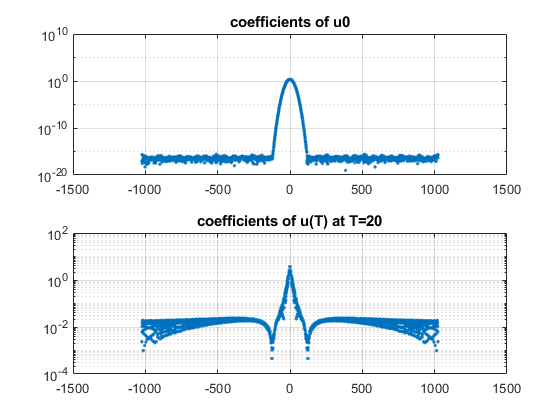}
\includegraphics[width=0.49\textwidth]{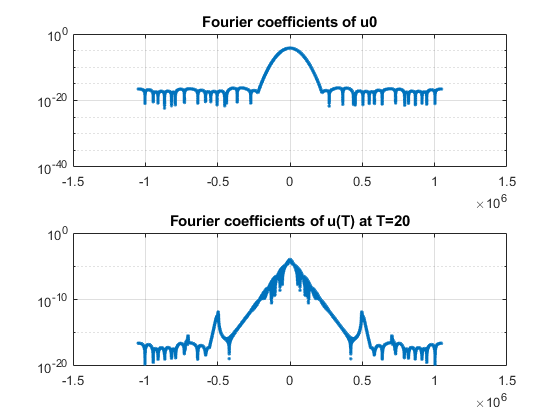}
\caption{Case $m=2$ (BO) with the initial condition $u_0=4e^{-x^2}$. Left: coefficients $a_n(t)$ in the rational basis expansion \eqref{rational}. Right: coefficients $a_n(t)$ in the Fourier basis expansion \eqref{Fourier}. Top row at time $t=0$, bottom row at $t=20$. One can note that the Fourier spectral method performs slightly better on the truncation error.}
\label{u m2 coefficients}
\end{figure}

Both Fourier and rational basis expansions work well for tracking solitons, however, the Fourier spectral method owns a better decay of the coefficients when tracking the scattering solutions for long time, see Figure \ref{u m2 coefficients}. The coefficients of $u(t)$ at $t=T$ (here $T=20$) for rational basis expansion remain at $10^{-2}$ beyond $N=1000$, which may lead to a large numerical error if not enough terms are used. This result is similar to \cite{TM98}. One possible reason for this is that the solution suffers from the under-resolution issue by using the rational basis. Therefore, unless otherwise specified, the numerical results reported below are obtained from the Fourier discretization. Moreover, as mentioned before, the Fourier spectral method enables us to  apply the mETDRK4, which allows us to choose a larger time step $\Delta t$, and thus, saves the computational time significantly.

\subsubsection{Scattering in a large domain} \label{scattering}
For some cases, where scattering or radiation happens (for example, for the initial data $u_0=\frac{1}{1+x^p}$ with $p=2$ or $4$ in the BO equation), a very large simulation domain is needed to track solutions behavior. Consequently, the rational basis for the space discretization is preferred.

\begin{figure}[ht]
\includegraphics[width=0.32\textwidth]{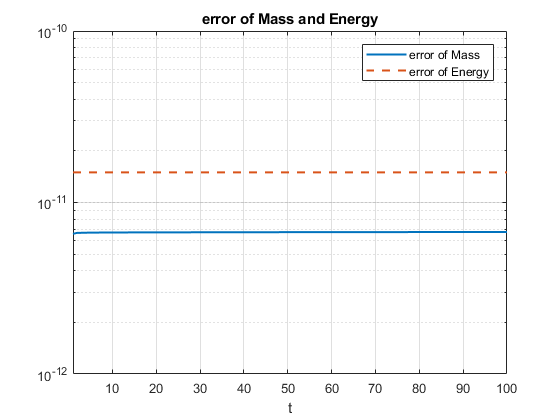}
\includegraphics[width=0.32\textwidth]{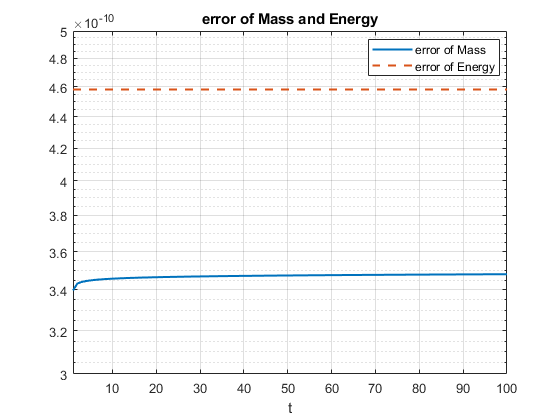}
\includegraphics[width=0.32\textwidth]{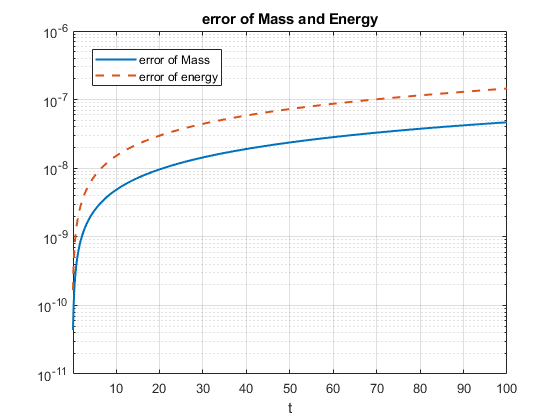}
\caption{Error of Mass and Energy \eqref{ME-error} in the case $m=2$ (BO) with $u_0=\frac{1}{1+x^2}$.  Left: from the Fourier spectral method with $N=2^{22}$, $L=20000\pi$, $\Delta t =5\times 10^{-3}$. Middle: from the Fourier spectral method with $N=2^{22}$, $L=20000\pi$, $\Delta t = 10^{-2}$. Right: from the rational basis with $N=2^{14}$, $L=800$, $\Delta t=4\times 10^{-3}$.}
\label{m2 ee em}
\end{figure}

To check the consistency, we track the error of mass and energy with respect to time
\begin{equation}\label{ME-error}
\mathcal{M}^{(k)}=\max M[u^{(k)}]-\min M[u^{(k)}] \quad \mbox{and} \quad \mathcal{E}^{(k)}=\max E[u^{(k)}]-\min E[u^{(k)}].
\end{equation}
We show these errors for the initial data $u_0 = \frac1{1+x^2}$ evolving with the BO time evolution in Figure \ref{m2 ee em} comparing both Fouier and rational basis expansions.

When using the Fourier spectral method in space and ETDRK4 in time, we set $N=2^{22}$, $L=20000 \pi$, and show that the errors of the mass and energy, $\mathcal{M}^{(k)}$ and $\mathcal{E}^{(k)}$, respectively, stay on the order of $ 10^{-11}$ when taking $5\times \Delta t=10^{-3}$, and the quantities $\mathcal{M}^{(k)}$ and $\mathcal{E}^{(k)}$ stays on the order $ 10^{-10}$ if taking $\Delta t=10^{-2}$ (see left and middle plots in Figure \ref{m2 ee em}).

When applying the rational basis functions in space and RK4 in time, due to the CFL condition, when taking the number of nodes $N=2^{14}$ and the mapping parameter $L=800$, the time step has to be on the order of $\Delta t=4 \times 10^{-3}$ to keep the numerical scheme in the stability regime. With the same initial condition $u_0=\frac{1}{1+x^2}$, we obtain the quantities $\mathcal{M}^{(k)}$ and $\mathcal{E}^{(k)}$ stabilized on the order of $10^{-7}$ (see the right plot in Figure \ref{m2 ee em}). This error does not decrease if we shrink the time step size $\Delta t$. However, it decreases if we increase the number of nodes $N$, consequently, we also have to shrink the time step size $\Delta t$. This justifies that the error is mainly caused by the under-resolution issue.

\subsection{The BO equation: toward the soliton resolution} \label{S:BO}
In this section, we report our results for the numerical simulations in the BO case $m=2$.
We first recall the rescaled ground state $Q_c(x)=\frac{4c}{1+(cx)^2}$, $c>0$, which will be important in the soliton shape fitting. Note that the height of $Q_c$ is $4$ times the constant $c$, so we will be measuring the height of the solution divided by $4$.

\begin{figure}[ht]
\includegraphics[width=0.49\textwidth]{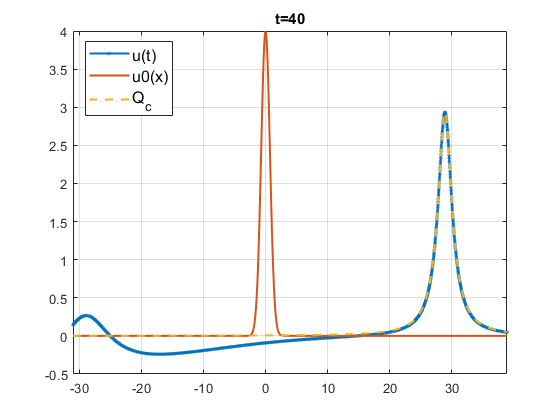}
\includegraphics[width=0.49\textwidth]{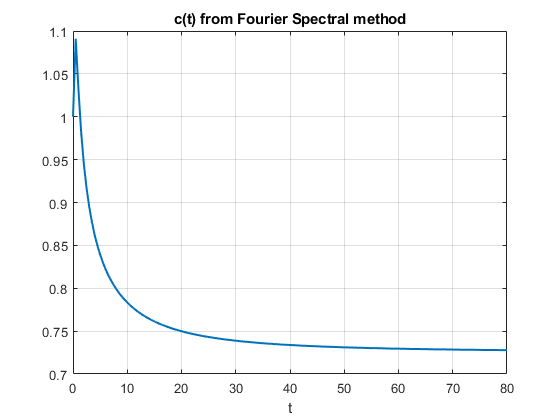}
\includegraphics[width=0.49\textwidth]{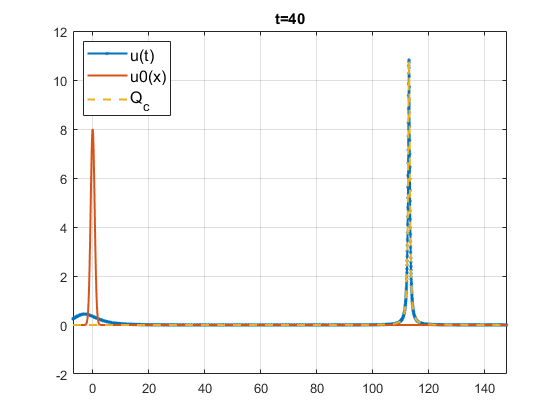}
\includegraphics[width=0.49\textwidth]{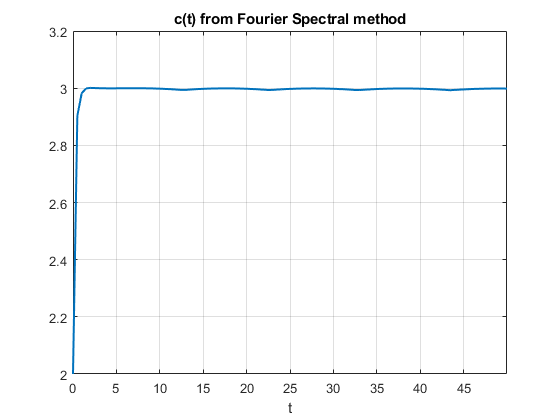}
\caption{Case $m=2$ (BO), evolution of gaussian initial data. Left top: $u_0=4e^{-x^2}$ (red) evolves (blue) into a soliton moving to the right (fitted in orange) at time $t=40$. Right top: decrease and leveling off of $c(t) = \|u(t)\|_{L^\infty}/4$ until $t=80$.  Left bottom: $u_0=8e^{-x^2}$ (red) evolves (blue) into a (faster) moving soliton (orange) at time $t=40$. Right bottom: $c(t)$ rapid growth and leveling off, $0<t<50$.}
\label{m2 Gaussian A8}
\end{figure}

\begin{figure}[ht]
\includegraphics[width=0.49\textwidth]{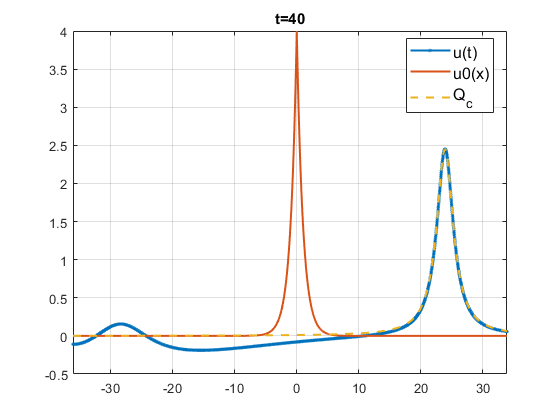}
\includegraphics[width=0.49\textwidth]{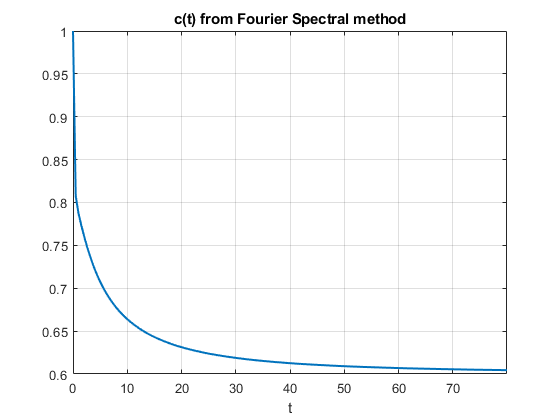}
\includegraphics[width=0.49\textwidth]{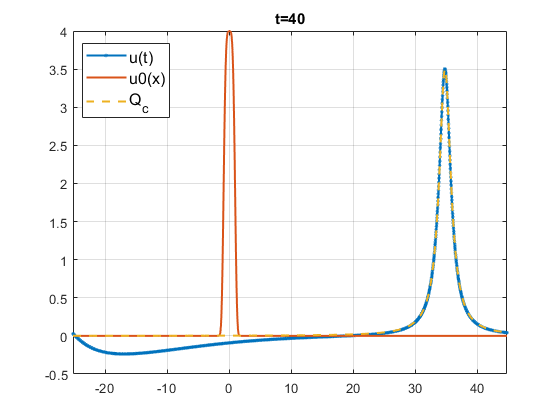}
\includegraphics[width=0.49\textwidth]{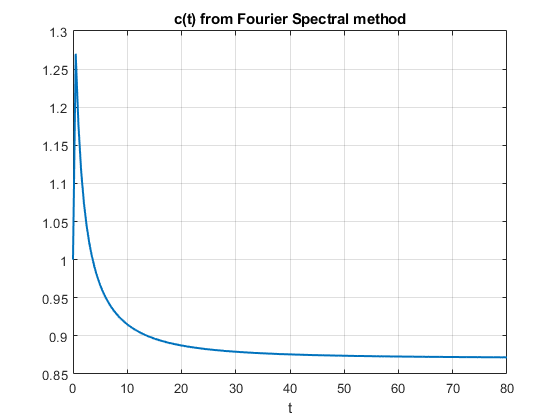}
\caption{Case $m=2$ (BO), evolution of exponentially decaying initial data. Left top: $u_0=4e^{-|x|}$ (red) evolves (blue) into a (small, slower) soliton moving to the right (fitted in orange) at time $t=40$. Right top: decrease and leveling off of $c(t)$ until $t=80$.   Left bottom: superGaussian $u_0=4e^{-x^4}$ (red) evolving (blue) into a soliton (orange) moving to the right at time $t=40$. Right: decrease and leveling off of $c(t)$ until $t=80$.}
\label{m2 supper-Gaussian}
\end{figure}

\begin{figure}[ht]
\includegraphics[width=0.49\textwidth]{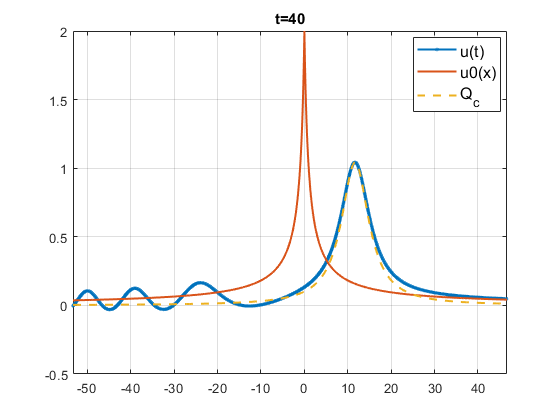}
\includegraphics[width=0.49\textwidth]{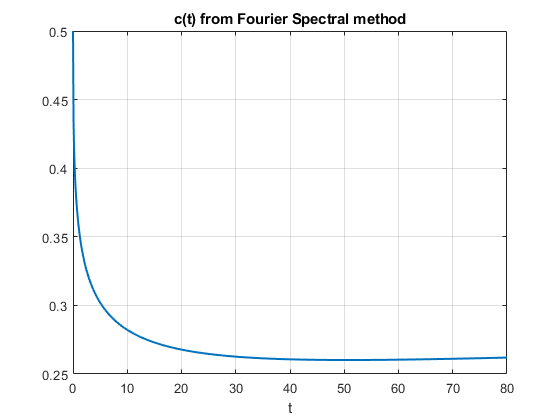}
\includegraphics[width=0.49\textwidth]{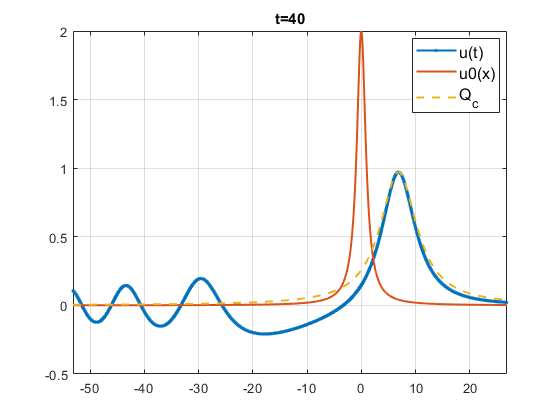}
\includegraphics[width=0.49\textwidth]{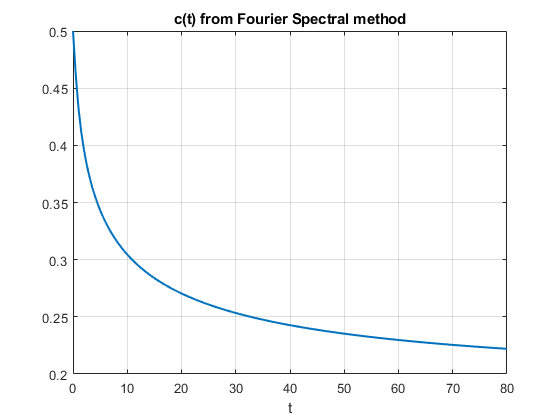}
\includegraphics[width=0.49\textwidth]{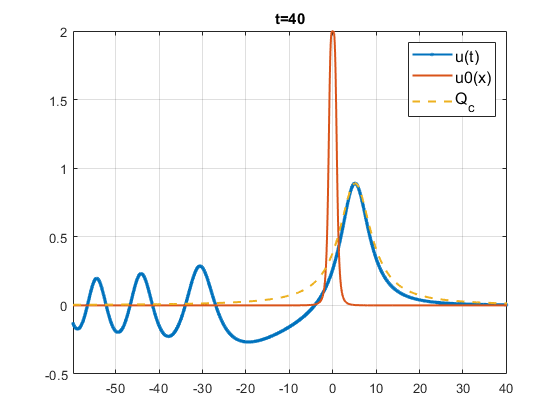}
\includegraphics[width=0.49\textwidth]{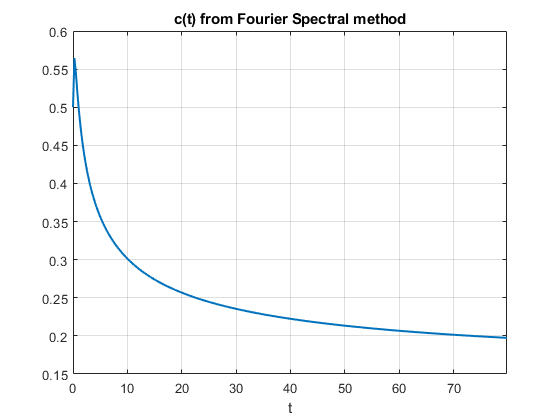}
\caption{Case $m=2$ (BO), evolution of polynomially decaying data.
Left top: $u_0=\frac{2}{1+|x|}$ (red) evolves (blue) into a soliton (fitted in orange) at time $t=40$. 
Left middle: $u_0=\frac{2}{1+x^2}$ (red) evolves (blue) into a (slower) soliton (orange) at time $t=40$.
Left bottom: $u_0=\frac{2}{1+x^4}$ (red) evolves into a (slow) soliton (orange) at time $t=40$. All right plots: decrease and leveling off of $c(t)$, respectively to each case, until $t=80$.}
\label{m2 poly-x2}
\end{figure}

\begin{figure}[ht]
\includegraphics[width=0.47\textwidth]{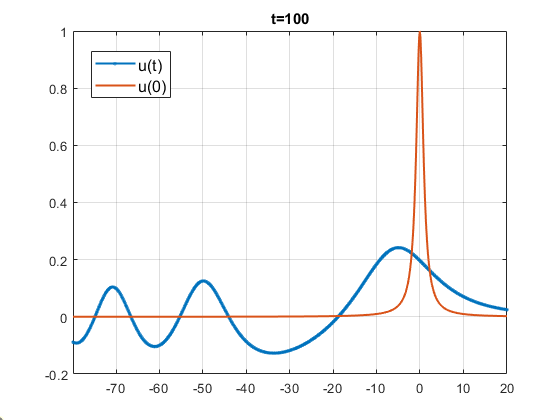}
\includegraphics[width=0.47\textwidth]{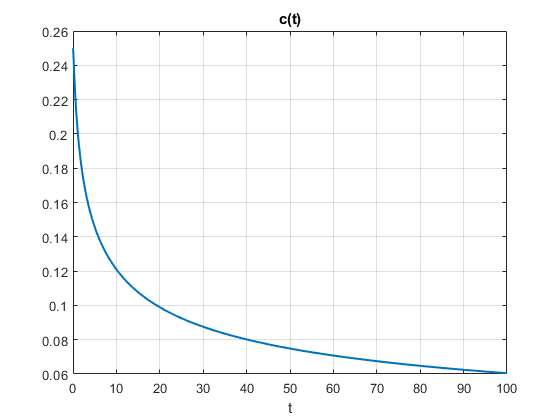}
\includegraphics[width=0.47\textwidth]{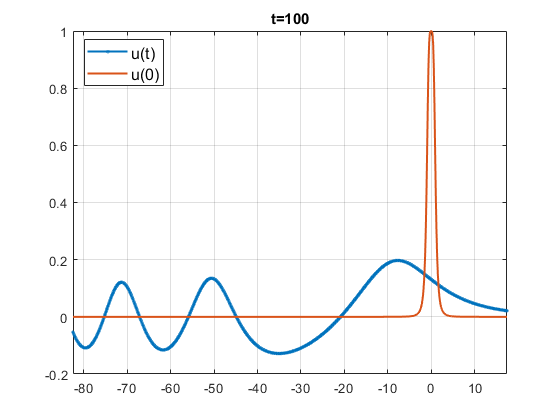}
\includegraphics[width=0.47\textwidth]{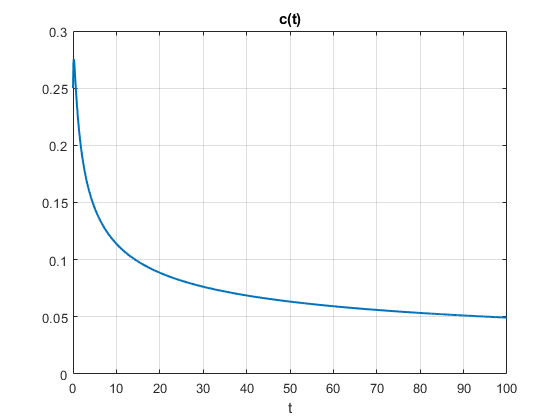}
\caption{Case $m=2$ (BO), evolution of polynomial decay data with small amplitude,
$u_0={1}/({1+x^2})$ (top) and $u_0=
{1}/({1+x^4})$ (bottom) (red) evolving into  radiation (blue) going to the left, at $t=100$. 
Right: decay of $c(t)$ up to time $t=100$. }
\label{m2 scatter}
\includegraphics[width=0.47\textwidth]{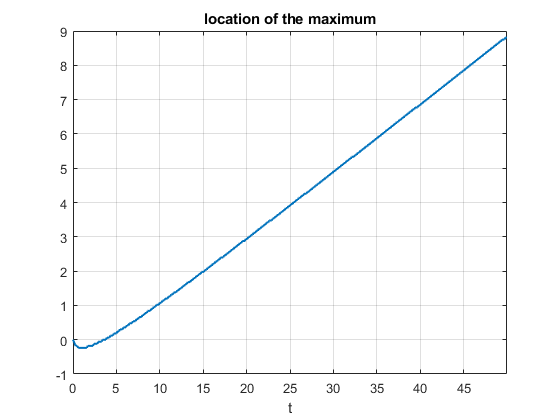}
\includegraphics[width=0.47\textwidth]{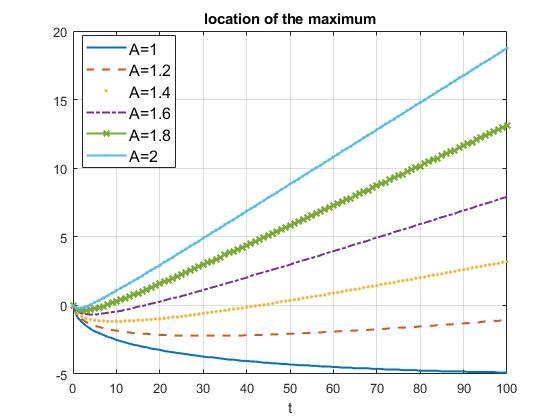}
\caption{Case $m=2$ (BO). Left: the evolution of $x_c(t)$ for $u_0=2/(1+x^2)$. 
Right: $x_c(t)$ trajectories for different A for $u_0=A/(1+x^2)$. 
The linear growth indicates the constant speed of the rescaled soliton formed from the initial data.}
\label{m2 poly-dichotomy}
\end{figure}

By taking various types of initial data, we study how their BO flow resolves into a soliton moving to the right and radiation going to the left.
Figures \ref{m2 Gaussian A8}--\ref{m2 scatter} show the time evolution of different initial data, which either approaches a rescaled and shifted soliton moving to the right and radiation going to the left, or just the radiation starting to go to the left (see further comments on this later). Each figure is split into two parts: (i) the left plot shows the solution profile: the initial condition is in red, its time evolution at the specific given time is in blue, and the fitting with the rescaled and shifted soliton profile $Q_c$ is in orange; (ii) the right plot shows the evolution of the magnitude, or the rescaling parameter, or soliton speed, $c(t)=\|u(x)\|_{L^\infty}/4$ as time increases.

\begin{figure}[ht]
\includegraphics[width=0.49\textwidth]{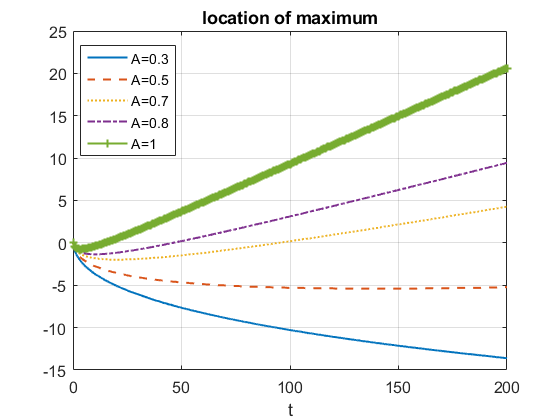}
\includegraphics[width=0.49\textwidth]{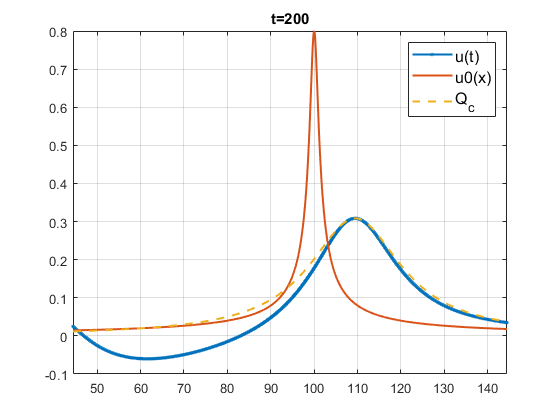}
\caption{Case $m=2$ (BO). Left: the center $x_c(t)$ trajectories for different $A$ for $u_0=A/\sqrt{1+x^2}$ up to $t=200$. The linear growth indicates the constant speed of the rescaled soliton formed from the initial data. Right: the solution profile $u(x,t)$ (blue) from
$u_0=0.8/\sqrt{1+x^2}$ (red) at $t=200$ with the emerging soliton profile $Q_c$ (orange).
}
\label{m2 poly-dichotomy2}
\end{figure}

We study the time evolution of the initial data with Gaussian decay in Figure \ref{m2 Gaussian A8}, with the general exponential decay
 in Figure \ref{m2 supper-Gaussian}, with polynomial decay but sufficiently large magnitude in Figure \ref{m2 poly-x2}, and  with polynomial decay but with smaller amplitude in Figure \ref{m2 scatter}. We consider a single maximum initial conditions, monotonously decaying to infinity, and we observe that if the solution (more precisely, the center or the location of the maximum) starts moving to the right, it will continue to do so. That indicates that it will resolve into a soliton\footnote{It is possible that it may resolve into several solitons; for smooth rapidly decaying data it should be possible to check via the IST by examining the discrete spectrum of the corresponding potential.} going to the right (or asymptotically approach a rescaled and shifted soliton) and radiation going to the left. If a center starts moving to the left, than one can say that the solution radiates to the left, however, we observe interesting dynamics that we discuss later and show in Figures \ref{m2 poly-dichotomy} and \ref{m2 poly-dichotomy2}.

If $c(t)$ converges to a horizontal asymptote as $t \rightarrow \infty$, this indicates that the solution $u(x,t)$ finally converges to a shifted version of the rescaled soliton $Q_c$, as was proved in \cite{KM2009}. On the other hand, if the quantity $c(t)$ keeps decreasing with respect to time $t$, this indicates that the solution scatters to zero. The right plots in Figures \ref{m2 Gaussian A8} -- \ref{m2 poly-x2} show that $c(t)$ converges to a horizontal asymptote, thus, showing the asymptotic stability of solitons and supporting the soliton resolution conjecture.
Additionally, by comparing the left subplots in Figure \ref{m2 Gaussian A8}, we can deduce that the  larger amplitude initial data leads to a faster convergence to the soliton $Q_c$, than the smaller amplitude amplitude initial data.

We next take smaller initial data and with polynomial decay, $u_0 = \frac1{1+x^p}$, $p=2,4$. In Figure \ref{m2 scatter} (top left), we observe that the evolution $u(x,t)$ from $u_0=\frac{1}{1+x^2}$ is starting to radiate to the left (note that the center of the maximum is traveling to the left). The right top plot there shows that $c(t)$ keeps decreasing even for an extended computational time up to $t=100$. 
A similar behavior happens when we take $u_0=\frac{1}{1+x^4}$, see Figure \ref{m2 scatter} (bottom). One could make a hypothesis that these solutions would finally all radiate to the left, or in other words, scatter to zero.

A natural question then is what exactly happens to the time evolution when a single peak data starts moving to the left. To investigate this numerically, we denote by $x_c$  the $x$-coordinate of the peak (or maximal amplitude) of the solution $u(x,t)$ at time $t$, i.e., $u(x_c,t)= \| u(x,t) \|_{L^{\infty}_x}$ and study the trajectory of $x_c$.
For that, we fix the decay of the initial data, for example, as $|x|^{-2}$, and track the trajectory of $x_c$. We start with $u_0= \frac2{1+x^2}$. The left plot in Figure \ref{m2 poly-dichotomy} shows that the value of $x_c$ grows linearly with respect to time $t$ (though notice a small decrease to the negative values initially, thus, indicating that the solution exhibits the radiative behavior initially). However, after a short time, the solution resolves into a soliton (plus radiation) and moves with a constant speed, matching the predicted soliton solution $Q_c(x-ct + C)=\frac{4c}{1+c^2(x-ct +C)^2}$.

Next, we take the same initial decay but vary the amplitude: $u_0(x) = \frac{A}{1+x^2}$, $A > 0$. With that data we track the $x_c(t)$ for different value of $A$. In the right plot of Figure \ref{m2 poly-dichotomy}, one can observe that a solution, initially traveling to the left, will reverse the direction of travel of $x_c$ and start going to the right, emerging a rescaled soliton (for example, as seen for $A=1.2$ or $1.4$). We do not see such behavior for $A=1$, simply because our computational time in that example was $0<t<100$. To investigate that further, we take $u_0 = \frac{1}{\sqrt{1+x^2}}$, a slower decay profile and simulate its time evolution for various values of the parameter $A$, see Figure \ref{m2 poly-dichotomy2} left plot. There we were able to run the simulations up to time $t=200$. We notice a similar behavior of the peak location for $A > 0.3$ (for smaller $A$ it is inconclusive due to the insufficient computational time; for slow decaying data IST methods are not available either). If we consider the evolution of the initial condition with $A=0.8$ (right plot in Figure \ref{m2 poly-dichotomy2}), we observe that indeed a soliton has emerged  - the fitting of the solution and the soliton now on the right is apparent, again confirming the soliton resolution.

\subsection{The mBO equation:  near the threshold behavior} \label{S:mBO}
When $m=3$, the equation \eqref{gBO} is $L^2$-critical, and it is
expected to have a stable blow-up.
In this part, as we do not anticipate simulations run for a very long time (due to blow-up), when applying the rational basis in space and the standard RK4 in time, we set the number of nodes $N=2048$, the mapping parameter $L=20$ and the time step $\Delta t=10^{-4}$; when applying the Fourier spectral method in space and the mETDRK4 in time, we set the nodes $N=2^{22}$, domain length $L=4000 \pi$ and the time step $\Delta t=10^{-3}$. We remark that we obtain very similar results from these two different approaches.

\begin{figure}[ht]
\includegraphics[width=0.49\textwidth]{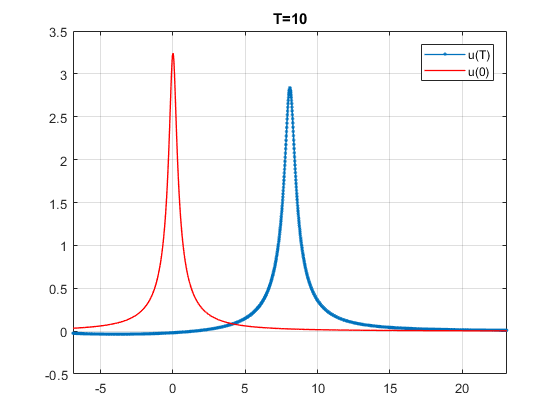}
\includegraphics[width=0.49\textwidth]{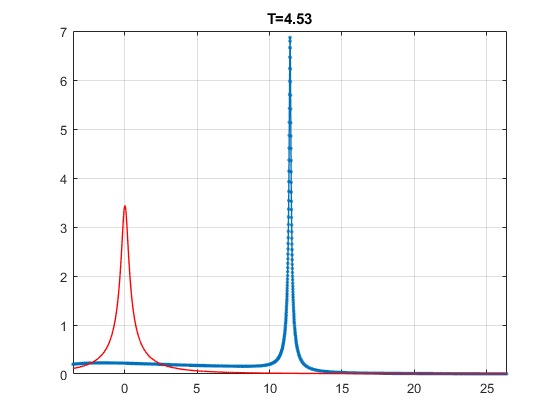}
\includegraphics[width=0.49\textwidth]{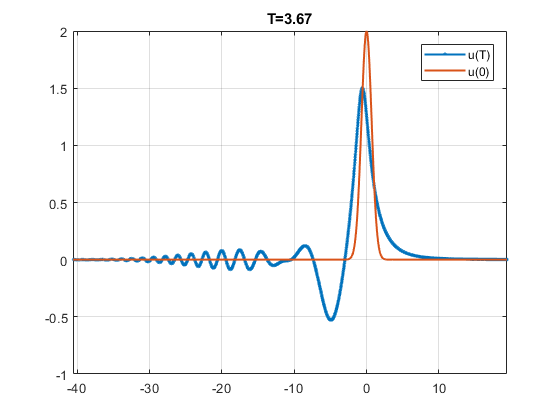}
\includegraphics[width=0.49\textwidth]{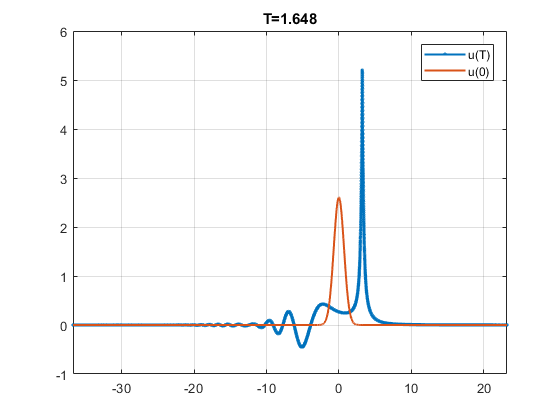}
\caption{Case $m=3$ (mBO). Left: scattering from $u_0=0.99Q$ (top) and $u_0=2e^{-x^2}$(bottom). Right: blow-up solutions $u_0=1.05Q$ (top) and $u_0=2.7e^{-x^2}$(bottom).}
\label{m3 profile}
\end{figure}

We start with considering the data $u_0 = A \, Q$, with $A \approx 1$ (for example, since it was shown in \cite{MP2016} that $Q$ is unstable).
Figure \ref{m3 profile} shows a scattering solution with the subcritical initial mass $u_0=0.99Q$ (top left) and a blow-up solution with supercritical initial mass $u_0=1.05Q$ (top right). The bottom row of that Figure indicates scattering and blow-up for the exponentially decaying data: the subcritical mass initial data, e.g., $u_0=2e^{-x^2}$, leads to scattering, and the supercritical amount of mass in the initial data, e.g., $u_0=2.7e^{-x^2}$, will lead toward the blow up behavior, see Figure \ref{m3 profile} (bottom). This gives partial confirmation of  Conjecture \ref{C: dichotomy}.

\begin{figure}[ht]
\includegraphics[width=0.32\textwidth]{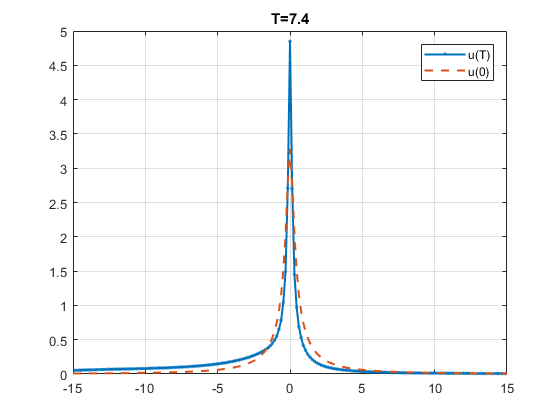}
\includegraphics[width=0.32\textwidth]{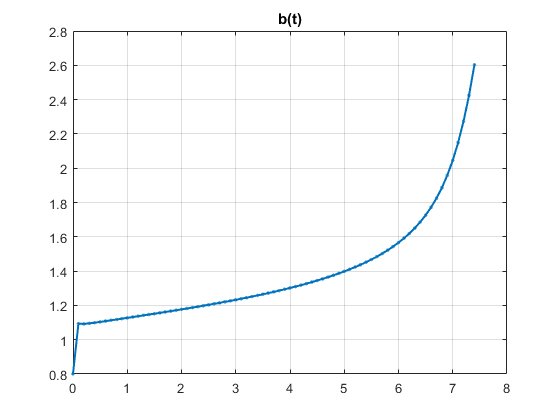}
\includegraphics[width=0.32\textwidth]{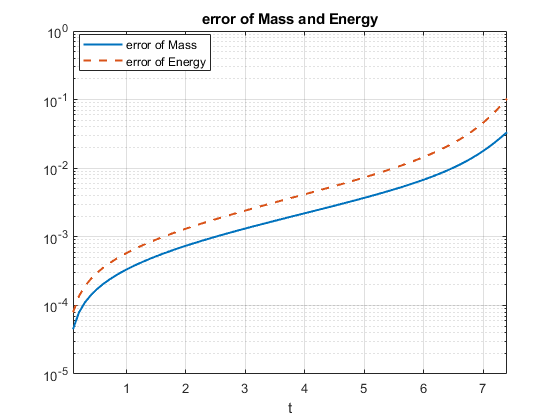}
\caption{Case $m=3$ (mBO). Left: profile of the solution for the re-centered equation with $u_0=1.01Q$. Middle: the speed of the solution, $b(t)$. Right: error of mass and energy.}
\label{m3 recenter}
\includegraphics[width=0.49\textwidth]{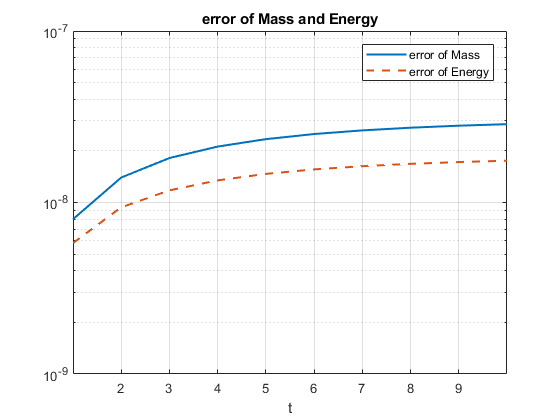}
\includegraphics[width=0.49\textwidth]{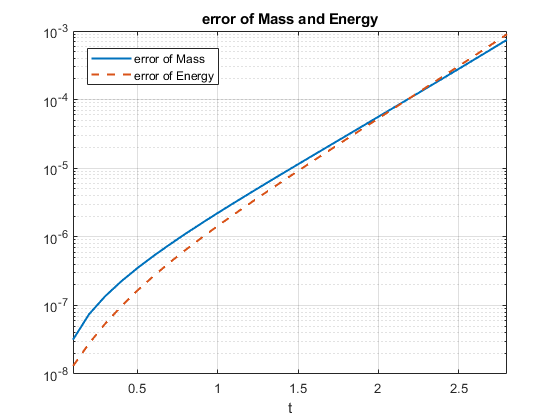}
\caption{Case $m=3$ (mBO). Left: error of mass and energy for scattering case. Right: error of mass energy for blow-up case. Initial data used: $u_0=0.99Q$ and $u_0=1.05Q$, respectively.}
\label{m3 em ee}
\end{figure}

It is worth noting that the blow-up behavior is numerically challenging to capture. For example, via direct simulations of solutions to the equation \eqref{BO}, the smallest supercritical mass that can be reliably observed to blow up is $u_0=1.05Q$. This is largely due to the under-resolution issue as well as the discretization error on the operator $\mathcal H \partial_{xx}$ (see similar issues and explanations for that in the NLS equation \cite{Fibich2015}, \cite{DM2002}). In order get closer to the threshold boundary in simulations, we re-center the equation (or use the ``commoving frame") by introducing a new variable $v(\xi,t)=u(x-x_c(t),t)$ such that $v$ stays at the center. Taking $\xi=x-x_c(t)$, the equation \eqref{BO} changes to
\begin{align}\label{BO_rescaled}
v_t-(\partial_tx_c)v_{\xi}+ \mathcal H v_{\xi \xi}+\frac{1}{m}\left( v^m \right)_{\xi}=0.
\end{align}
By denoting $b=\partial_t x_c$, one can observe that $b$ is the speed of the solution. Tracking the parameter $b$ shows how the speed changes during the time evolution.
On the other hand, to represent the quantity $b$ in terms of $v,\xi,t$, we adopt the argument from \cite{KP13,LSS13}. Since $x_c$ is the location of the maximum point of $u(x,t)$ for any time $t$ (or $v(\xi,t)$ at $\xi=0$, equivalently), one sufficient condition is that $v_{\xi}\big \vert_{\xi=0}=0 $ for all $t$. Differentiating \eqref{BO_rescaled} with respect to $\xi$ and evaluating at $\xi=0$, one has
\begin{align}\label{BO_rescaled b}
b=\dfrac{\mathcal H v_{\xi \xi \xi}+\frac{1}{m} (v^m)_{\xi \xi}}{v_{\xi \xi}}\Big \vert_{\xi=0}.
\end{align}
In our simulation, the equations \eqref{BO_rescaled} and \eqref{BO_rescaled b} are solved simultaneously by applying the rational basis functions in space discretization and standard RK4 \eqref{RK4} in time evolution. Using this reformulation, we can circumvent the under-resolution issue, and show that the solution blows up in finite time, for example, with $u_0=1.01Q$, see Figure \ref{m3 recenter}.

We also tested the Gaussian type initial data $u_0=Ae^{-x^2}$ and obtained similar results, as we have shown in Figure \ref{m3 profile} (bottom). Setting $u_0=Ae^{-x^2}$, the condition $M[u_0]<M[Q]$ is satisfied when $A<2.4312$ (and we observe that such data scatters). On the other hand, $E[u_0]<0$ when $A>2.6020$. We actually observe the blow-up when $A \geq A^*=2.55$, the energy for the given $A^*$, $E[u_0]$ is positive. These results numerically validate Conjecture \ref{C: dichotomy} (and remark afterward).

To finish this subsection, we mention that Figure \ref{m3 em ee} tracks the error of mass and energy for the scattering and blow-up cases, respectively. The scheme maintains accuracy for scattering even for the long time simulations, but may become less reliable when the solution blows up. It is mainly due to the numerical difficulties in dealing with blow-up solutions.  However, up to the time we end our simulations ($\sim T=4.5$), these error quantities remain at a reasonable level of accuracy for the blow-up behavior.

\subsection{The gBO equations: blow-up vs global existence dichotomy} \label{S:gBO}
\begin{figure}[ht]
\includegraphics[width=0.49\textwidth]{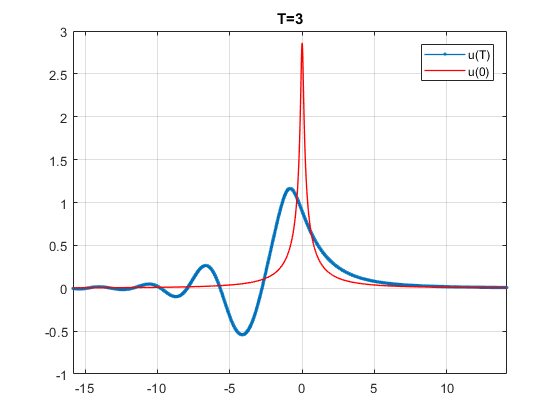}
\includegraphics[width=0.49\textwidth]{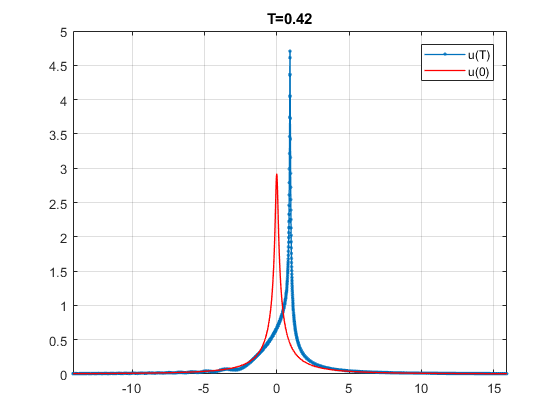}
\includegraphics[width=0.49\textwidth]{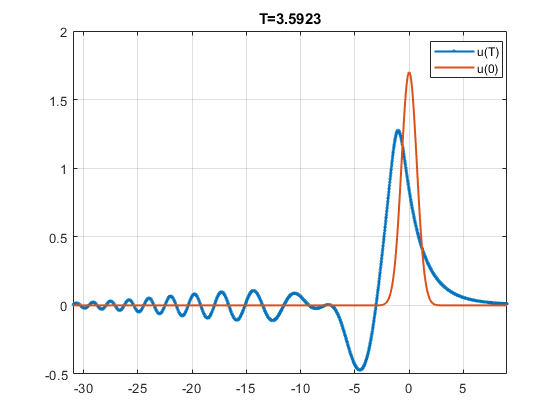}
\includegraphics[width=0.49\textwidth]{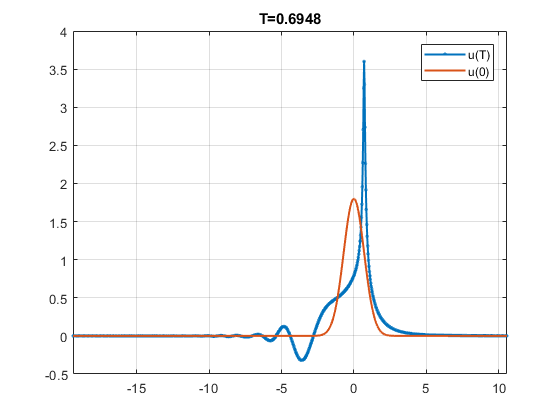}
\caption{Case $m=4$ (gBO). Left: radiating to the left $u_0=0.99Q$ (top), $u_0=1.7e^{-x^2}$ (bottom). Right: blow-up solutions $u_0=1.01Q$ (top), $u_0=1.8e^{-x^2}$ (bottom).}
\label{m4 profile}
\includegraphics[width=0.49\textwidth]{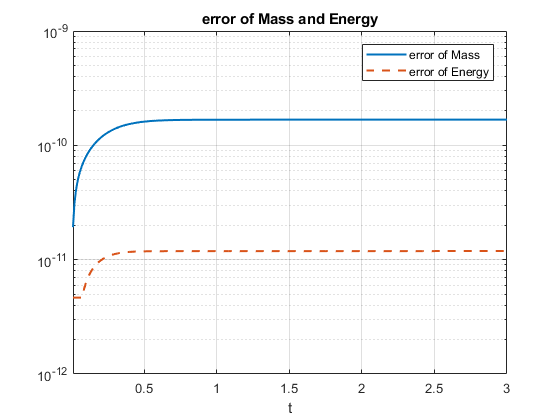}
\includegraphics[width=0.49\textwidth]{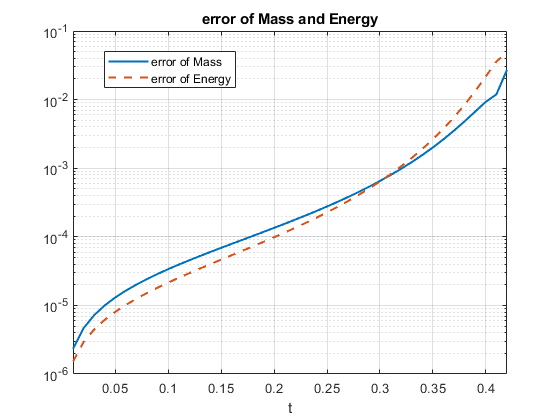}
\caption{Case $m=4$ (gBO). Error of mass and energy in the scattering case (left) and in the blow-up case (right); from $u_0=0.99Q$ and $u_0=1.01Q$, respectively.}
\label{m4 em ee}
\end{figure}

\begin{figure}[ht]
\includegraphics[width=0.49\textwidth]{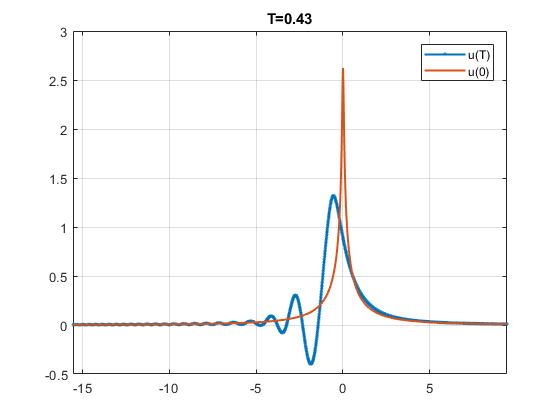}
\includegraphics[width=0.49\textwidth]{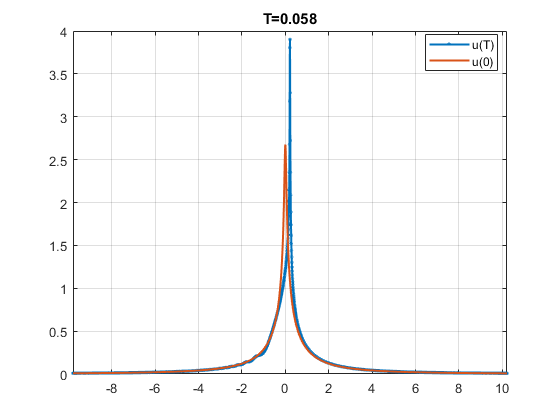}
\includegraphics[width=0.49\textwidth]{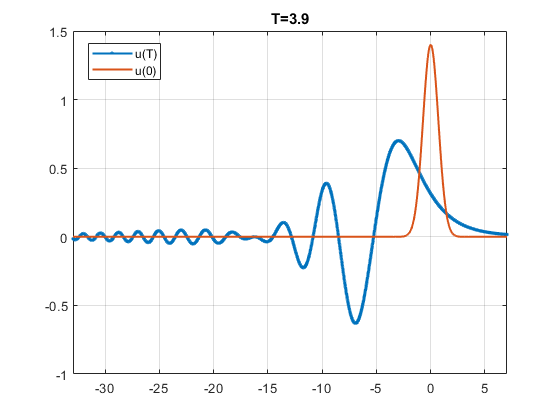}
\includegraphics[width=0.49\textwidth]{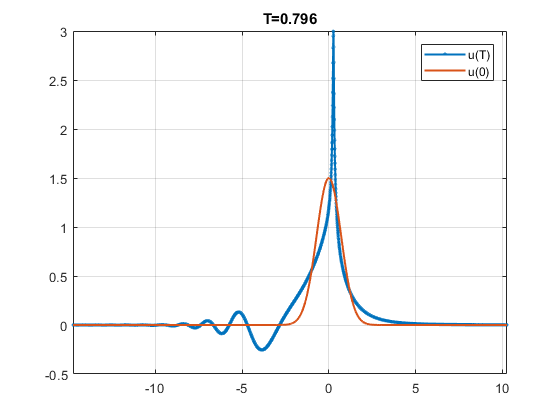}
\caption{Case $m=5$ (gBO). Left: radiating to the left $u_0=0.99Q$ (top), $u_0=1.4e^{-x^2}$ (bottom). Right: blow-up solutions $u_0=1.01Q$ (top), $u_0=1.5e^{-x^2}$ (bottom).}
\label{m5 profile}
\includegraphics[width=0.49\textwidth]{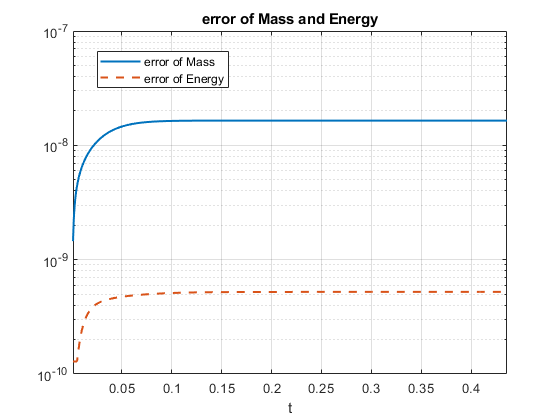}
\includegraphics[width=0.49\textwidth]{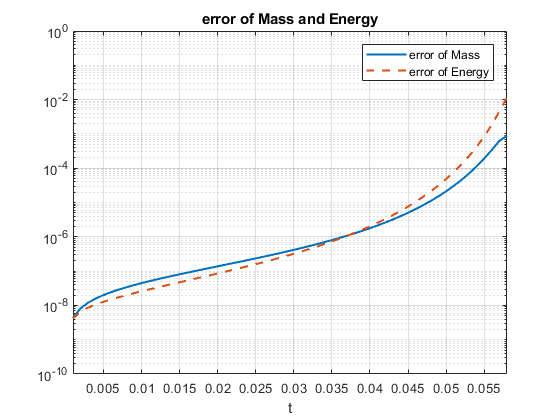}
\caption{Case $m=5$ (gBO). Error of mass and energy in the scattering case (left) and in the blow-up case (right); initial data is $u_0=0.99Q$ and $u_0=1.01Q$, respectively.}
\label{m5 em ee}
\end{figure}

Here we study solutions to the gBO equation \eqref{gBO} with $m=4,5$, the $L^2$-supercritical case. 
We consider a family of initial data $u_0 = A \, Q$, $A>0$, where $Q$ is the corresponding ground state solution of \eqref{Q} and plot the corresponding gBO time evolution with $A \approx 1$ in Figures \ref{m4 profile} and \ref{m5 profile} for $m=4$ and $5$, correspondingly. Besides that we also consider exponentially decaying data such as Gaussian, and indicate snapshots of time evolution in the above Figures.

In Figure \ref{m4 profile} ($m=4$) we observe that initial data such as $u_0=0.99Q$ and  $u_0=1.7e^{-x^2}$ start scattering, and initial data such as $u_0=1.01Q$ and $u_0=1.8e^{-x^2}$ shows the blow up behavior. Similar behavior can be seen in Figure \ref{m5 profile} for $m=5$. In Figures \ref{m4 em ee} and \ref{m5 em ee}, we track the errors of mass and energy in the corresponding cases, $m=4$ and $m=5$, using $u_0=0.99Q$ $u_0=1.01Q$ as initial conditions. Similar to the mBO case ($m=3$), the numerical scheme here maintains the same  accuracy for scattering solutions but becomes less reliable when solutions blow up. We terminate the simulation before the error becomes too large and mark the solution as blow-up for these initial conditions.

We turn to the justification of Conjecture \ref{C: dichotomy2}. First of all, we observe that any negative energy initial data that we used in the supercritical case (a single maximum, monotonously decreasing) produced a blow-up behavior, confirming the Part I of the Conjecture \ref{C: dichotomy2}. 
We next focus on the numerical validation of Part IIb, since the Part IIa of Conjecture \ref{C: dichotomy2} was proved in \cite{FLP2014}. 
We denote $u_0=A \, v$, where $A$ is a constant and the function $v$ is taken from the set
\begin{equation}\label{E:dataset}
\left\{ ~e^{-x^4}, ~e^{-x^2}, ~\dfrac1{1+x^4}, ~\dfrac1{1+x^2}, ~\dfrac1{\sqrt{1+x^2}} ~\right\}.
\end{equation}
We also set $\theta=(\frac12-s)/{s}$ (thus, $\theta=2$ if $m=4$ and $\theta=1$ if $m=5$)
 and denote
$$
\alpha_m = M[Q^{(m)}]^{\theta} E[Q^{(m)}], \quad \beta_m = M[Q^{(m)}]^{\theta} D[Q^{(m)}],
$$
where the index $m$ indicates the corresponding ground state from \eqref{Q}.
For brevity we also denote
$$
M_v=M[v], \quad D_v=D[v], \quad \mbox{and} \quad P_v=\frac{1}{m(m+1)} \| v \|_{L^{m+1}}^{m+1}.
$$
With the above notation, the conditions \eqref{ME bound} and \eqref{MD scatter} in Conjecture \ref{C: dichotomy2} are reduced to
\begin{align}\label{ME bound A}
A^{2\theta} M_v^\theta \left( \dfrac{1}{2} A^{2} D_v -A^{m+1} P_v \right) < \alpha_m,
\end{align}
\begin{align}\label{MD scatter A}
A^{2 \theta +2} M_v^{\theta} D_v <\beta_m.
\end{align}
When $m=4$, the equation \eqref{ME bound A} becomes
\begin{align}\label{ME bound m4}
\frac{1}{2} A^6 M_v^2 D_v - A^{9} M_v^2 P_v - \alpha_4 <0.
\end{align}
Denoting $B=A^3$ and incorporating the fact that $A>0$, the inequality \eqref{ME bound m4} yields
\begin{align}\label{ME bound m4 B}
\frac{1}{2} B^2 M_v D_v-B^{3} M_v P_v - \alpha_{4} <0, ~~A> 0.  
\end{align}
Similarly, when $m=5$, from \eqref{ME bound A}, one has
\begin{align}\label{ME bound m5}
\frac1{2} {A^4} M_v D_v-A^{8} M_v P_v - \alpha_5 <0.
\end{align}
With $B=A^4$ the inequality \eqref{ME bound m5} yields
\begin{align}\label{ME bound m5 B}
\frac1{2} {B} M_v D_v-B^{2} M_v P_v - \alpha_5 <0, ~~ A>0.
\end{align}
It is easy to get the roots in \eqref{ME bound m4 B} and \eqref{ME bound m5 B}, and obtain the threshold values of $A$, for which inequalities hold, consequently. For both cases $m=4$ and $m=5$, we have two values of $A$, denoting them as $A_0^+$ and $A_0^-$ (the larger and the smaller) such that \eqref{ME bound m4 B} and \eqref{ME bound m5 B} hold: $A<A_0^-$ and $A>A_0^+$. We denote the value of $A$ obtained from \eqref{MD scatter A} as $A_1$. In a similar manner, we obtain the value of $A$, which is the threshold for the negative energy $E[u_0]<0$, denoting it by $A_E$ ($A>A_E$).

After computing the quantities $A_0^-$, $A_0^+$, $A_1$ and $A_E$ for different initial data $u_0=A \, v$ (observe that we always have $A_0^- < A_1 < A_0^+ < A_E$), we then check the corresponding time evolutions. From our numerical simulations we identify the threshold for blow-up vs global existence behavior up to two decimal places of accuracy and denote it by $A_T$,
that is, for $A \geq A_T$ we observe blow-up and for $A < A_T$ we observe a globally existing solution (typically radiating to the left, however, this should be investigated further). We put all these values into Tables \ref{Table: super-Gaussian} - \ref{Table: poly x1}.

\begin{figure}[ht]
\includegraphics[width=0.9\textwidth,height=0.28\textheight]{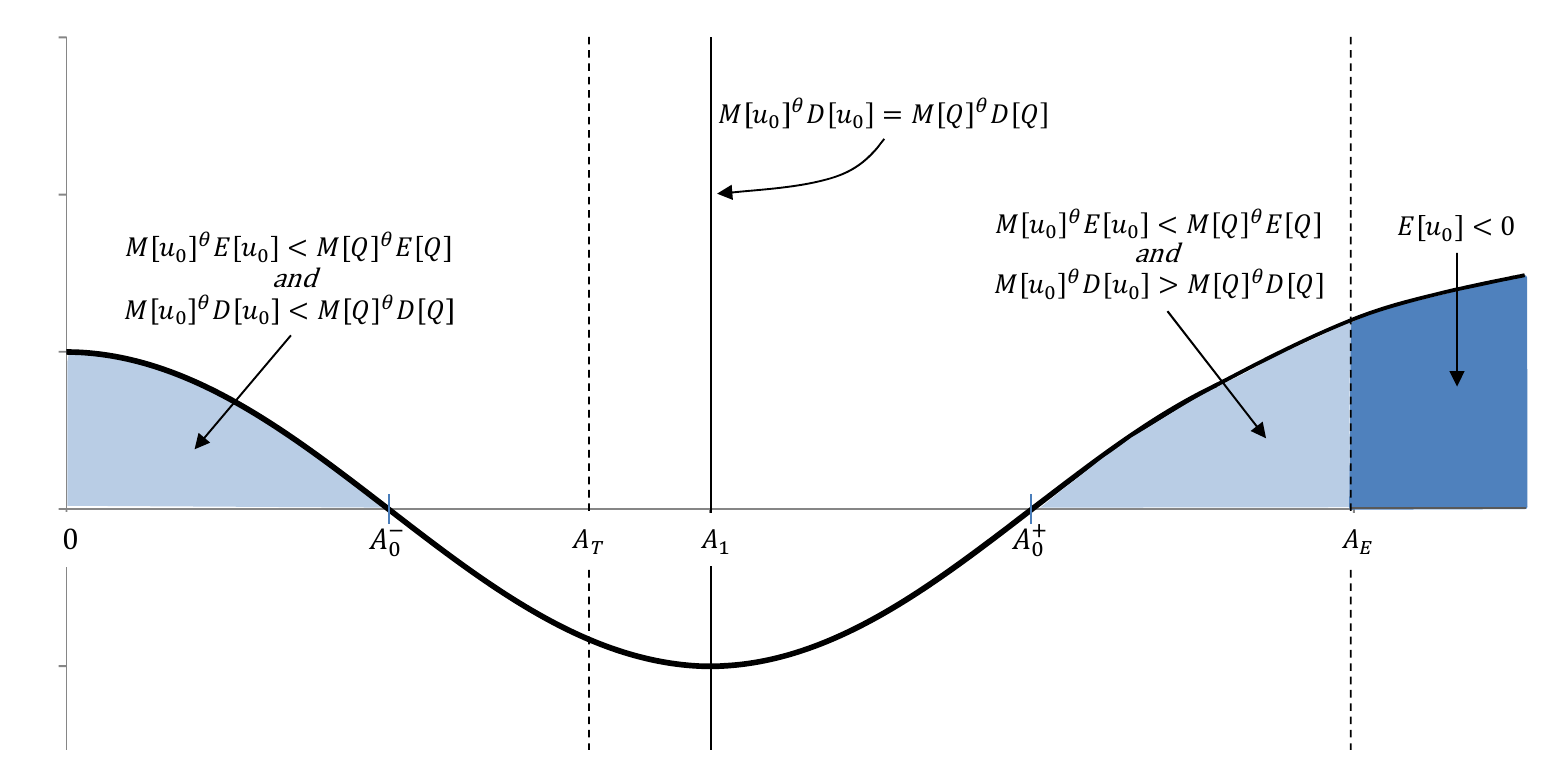}
\caption{Thresholds for $A$ values in the initial condition $u_0 = A \, v$ with $v$ as in \eqref{E:dataset} from Conjecture \ref{C: dichotomy2}, or \eqref{ME bound A}-\eqref{MD scatter A}, in the $L^2$-supercritical gBO equation.}
\label{fig_thresholds}
\end{figure}

To better understand the values in these tables, we provide a graphical representation of the computed thresholds for $A$ in Figure \ref{fig_thresholds}. This is a typical situation of the thresholds that we obtained, however, we remark that the values $A_T$ and $A_1$ can be interchanged for some data (e.g., $m=4$ in Table \ref{Table: poly x1}), nevertheless, it does not affect our conclusions. In all studied and computed cases, the obtained results that confirm Conjecture \ref{C: dichotomy2} (as well as it also supports Part IIa proved in \cite{FLP2014}).

\begin{table}
{\bf Threshold $A$ values for $u_0=Av$ for Conjecture \ref{C: dichotomy2}, or in \eqref{ME bound A}-\eqref{MD scatter A}, and numerical threshold $A_T$ for various functions $v$ in \eqref{E:dataset} in the $L^2$-supercritical gBO}
\\
\begin{tabular}{|l|l|l|l|l|l|}
\hline
$v=e^{-x^4}$&  $A_0^-$ & $A_0^+$ & $A_1$ & $A_E$ & $A_T$ \\
\hline
 $m=4$      & $1.4225$ & $2.2288$ & $1.6295$ & $2.2680$ & $1.57$ \\
 \hline
 $m=5$      & $1.1950$ & $2.0085$ & $1.3798$ & $2.0687$ & $1.32$ \\
 \hline
\end{tabular}
\caption{$v=e^{-x^4}$.}
\label{Table: super-Gaussian}
\begin{tabular}{|l|l|l|l|l|l|}
\hline
$v=e^{-x^2}$&  $A_0^-$ & $A_0^+$ & $A_1$ & $A_E$ & $A_T$ \\
\hline
 $m=4$      & $1.6518$  & $2.2462$ & $1.8428$ & $2.3279$ & $1.73$ \\
 \hline
 $m=5$      & $1.3982$ & $2.0278$ & $1.5801$ & $2.1338$ & $1.48$ \\
 \hline
\end{tabular}
\caption{$v=e^{-x^2}$.}
\label{Table: Gaussian}
\begin{tabular}{|l|l|l|l|l|l|}
\hline
$v=\frac{1}{1+x^4}$&  $A_0^-$ & $A_0^+$ & $A_1$ & $A_E$ & $A_T$ \\
\hline
 $m=4$      & $1.4495$ & $2.0494$ & $1.6308$ & $2.1107$ & $1.55$ \\
 \hline
 $m=5$      & $1.2405$ & $1.8795$ & $1.4125$ & $1.9630$ & $1.34$ \\
 \hline
\end{tabular}
\caption{$v=\frac{1}{1+x^4}$.}
\label{Table: poly x4}
\begin{tabular}{|l|l|l|l|l|l|}
\hline
$v=\frac{1}{1+x^2}$&  $A_0^-$ & $A_0^+$ & $A_1$ & $A_E$ & $A_T$ \\
\hline
 $m=4$      & $1.6668$ & $1.9446$ & $1.7795$ & $2.0910$ & $1.70$ \\
 \hline
 $m=5$      & $1.4556$ & $1.8107$ & $1.5864$ & $1.9758$ & $1.50$ \\
 \hline
\end{tabular}
\caption{$v=\frac{1}{1+x^2}$.}
\label{Table: poly x2}
\begin{tabular}{|l|l|l|l|l|l|}
\hline
$v=\frac{1}{(1+x^2)^{\frac{1}{2}}}$&  $A_0^-$ & $A_0^+$ & $A_1$ & $A_E$ & $A_T$ \\
\hline
 $m=4$      & $1.4151$ & $1.5195$ & $1.4627$ & $1.6839$ & $1.48$ \\
 \hline
 $m=5$      & $1.3461$ & $1.4819$ & $1.4059$ & $1.6873$ & $1.37$ \\
 \hline
\end{tabular}
\caption{
$v=\frac{1}{\sqrt{1+x^2}}$.}
\label{Table: poly x1}
\end{table}

\section{Conclusions}
In this work we considered the generalized Benjamin-Ono equation with different powers of nonlinearity. We first used the Petviashvili's iteration method to obtain the ground state solutions for different powers (coinciding with the explicit solution in $m=2$ case). Compared with other methods (e.g., \cite{BX12}), which use a homotopy technique, the Petviashvili's iteration method allows the use of a robust initial guess.

We then investigated behavior of solutions for various types of initial data (single maximum with  monotone decay) in the standard Benjamin-Ono equation by tracking the time evolution and fitting when possible with the rescaled and shifted solitons. In particular, we observed that the solution traveling to the right will approach a rescaled (and shifted) version of the soliton, thus, traveling with the speed approaching the speed of that rescaled soliton. Some solutions will start traveling to the left as radiation, however, we observe that (in the allowed computational time) the location of the maximum of such solution initially traveling to the left will slow down, then change the direction and start traveling to the right forming a rescaled version of the soliton. We were able to observe such behavior in exponentially decaying data and also in polynomially decaying data with the decay rate as slow as $1/|x|$ (this is a constrain of the numerical method with rational functions). This not only confirms the soliton resolution but also shows the intermediate process of soliton formation (and consequently, asymptotic stability of solitons in the BO equation).

We then studied the $L^2$-critical mBO equation, and numerically confirmed the existence of stable blow-up. In particular, we investigated the behavior of the initial data such as a multiple of the ground state and observed that, indeed, the mass of the ground state gives a threshold for the global existence vs blow-up in finite time in this case. We also observed that negative energy solutions (of different types of considered initial data) blow up in finite time.

Finally, we investigated the $L^2$-supercritical BO equation. A blow-up occurs for all negative energy (single maximum) initial data in our numerical simulations. Furthermore, we gave numerical confirmation for the blow-up vs global existence dichotomy under the mass-energy threshold (this also includes positive energy initial data). The results are similar to other  $L^2$-supercritical (or intercritical) cases of dispersive PDEs.

\bibliographystyle{abbrv}
\bibliography{ref}
\end{document}